\documentclass[12pt]{article}
\usepackage{amsmath,amsfonts,amssymb}

\usepackage{xcolor}
\usepackage{graphicx}

\usepackage[margin=1in]{geometry}
\usepackage[square,sort,comma,numbers]{natbib}

\usepackage{rotating}
\usepackage[linesnumbered,ruled]{algorithm2e}
\usepackage{caption}

\usepackage{tikz}
\usetikzlibrary{arrows,shapes}

\makeatletter
\pgfdeclareshape{datastore}{
  \inheritsavedanchors[from=rectangle]
  \inheritanchorborder[from=rectangle]
  \inheritanchor[from=rectangle]{center}
  \inheritanchor[from=rectangle]{base}
  \inheritanchor[from=rectangle]{north}
  \inheritanchor[from=rectangle]{north east}
  \inheritanchor[from=rectangle]{east}
  \inheritanchor[from=rectangle]{south east}
  \inheritanchor[from=rectangle]{south}
  \inheritanchor[from=rectangle]{south west}
  \inheritanchor[from=rectangle]{west}
  \inheritanchor[from=rectangle]{north west}
  \backgroundpath{
    \southwest \pgf@xa=\pgf@x \pgf@ya=\pgf@y
    \northeast \pgf@xb=\pgf@x \pgf@yb=\pgf@y
    \pgfpathmoveto{\pgfpoint{\pgf@xa}{\pgf@ya}}
    \pgfpathlineto{\pgfpoint{\pgf@xb}{\pgf@ya}}
    \pgfpathmoveto{\pgfpoint{\pgf@xb}{\pgf@yb}}
    \pgfpathlineto{\pgfpoint{\pgf@xb}{\pgf@ya}}
    \pgfpathmoveto{\pgfpoint{\pgf@xa}{\pgf@yb}}
    \pgfpathlineto{\pgfpoint{\pgf@xb}{\pgf@yb}}
 }
}

\newlength{\defbaselineskip}
\newcommand{\setlinespacing}[1]%
           {\setlength{\baselineskip}{#1 \defbaselineskip}}

\providecommand{\keywords}[1]
{
  \small	
  \noindent \textbf{Keywords:} \textit{#1}
}

\newtheorem{theorem}{Theorem}

\newtheorem{remark}{Remark}

\newcommand{\RR}{{\mathbb{R}}}

\newcommand{\sU}{{\cal U}}

\newcommand{\eqn}[1]{(\ref{#1})}
\newcommand{\beq}{\begin{displaymath}}
\newcommand{\eeqno}{\end{displaymath}}

\newcommand{\Var}{\hbox{Var}}

\newcommand{\qandq}{\quad\mbox{and}\quad}

\newcommand{\qforallq}{\quad\mbox{for all}\quad}

\newcommand{\beql}[1]{\begin{equation}\label{#1}}
\newcommand{\eeq}{\end{equation}}

\title{A Robust Queueing Network Analyzer \\ Based on  Indices of Dispersion}
\author{Ward Whitt\thanks{\scriptsize Department of Industrial Engineering and Operations Research, Columbia University}\\ \small ww2040@columbia.edu
\and
Wei You\thanks{\scriptsize Department of Industrial Engineering and Decision Analytics, Hong Kong University of Science and Technology}\\  \small weiyou@ust.hk}
\date{\today}

\begin{document}

\maketitle

\begin{abstract}

We develop a robust queueing network analyzer algorithm to approximate the steady-state performance of a single-class open queueing network of single-server queues with Markovian routing.
The algorithm allows non-renewal external arrival processes, general service-time distributions and customer feedback.
We focus on the customer flows, defined as the continuous-time processes counting customers flowing into or out of the network, or flowing from one queue to another.
Each flow is partially characterized by its rate and a continuous function that measures the stochastic variability over time.
This function is a scaled version of the variance-time curve, called the index of dispersion for counts (IDC).
The required IDC functions for the flows can be calculated from the model primitives, estimated from data or approximated by solving a set of linear equations.
A robust queueing technique is used to generate approximations of the mean steady-state performance at each queue from the IDC of the total arrival flow and the service specification at that queue.
The algorithm effectiveness is supported by extensive simulation studies and heavy-traffic limits.

\end{abstract}

\newpage

\section{Introduction}
This paper contributes to analytical methods for designing and optimizing service systems.
Such systems appear in a broad and diverse range of settings, including customer contact centers, hospitals, airlines, online marketplaces, ride-sharing platforms and cloud computing networks.
The design and operation of these systems is challenging,
largely because there is uncertainty about when customers will arrive and their service requirements.

Fortunately, useful guidance can often be provided by exploiting mathematical models using stochastic processes.
Prominent among these are stochastic queueing network models, because service is often provided in a sequence of steps;
e.g., see \cite{CY01,BvD11}.
There is an extensive literature on the applications of queueing network models to service systems.
For example, see \cite{SC81} for a review of applications in computer networks, see \cite{BJR15,FHS17,OW17} for examples in ride-sharing economies and see \cite{CL11,CDG16,ZA16,DS19,KWC18} for healthcare-related applications.

Service operation policies often rely on quantitative descriptions of the system performance,
called \textit{performance measures}, such as
the waiting time, the queue length, and the workload in the system.
Decision support for service operations relies on accurate characterization of these performance measures.

A standard way to analyze the performance of complex queueing models is to employ computer simulation, e.g., see \cite{SM05,Z11}.
However, as noted in \cite{DGS16}, a great disadvantage of simulation-based optimization methods is the often prohibitive computation time required to obtain optimal solutions for service operation problems involving a multidimensional stochastic network.  Thus,
analytical analysis of the models can be very helpful.
However, the class of queueing networks that can be solved analytically requires strong assumptions that are rarely satisfied, whereas more realistic models are prohibitively hard to analyze exactly.
Hence, analytical performance approximation of queueing networks remains an important tool.

In this paper, we provide convenient algorithms to approximate the steady-state performance measures in a single-class open queueing network (OQN) with Markovian routing, unlimited waiting space and the first-come first-served (FCFS) service discipline.
We focus on non-Markov OQNs where the external arrival processes need not be Poisson or renewal and the service-time
distributions need not be exponential.



\subsection{Approximation Algorithms}\label{sec:approx_algs}
In this section, we briefly review existing approximation algorithms for non-Markov OQNs;
additional literature review appears in an online appendix.

\subsubsection{Decomposition Approximations}\label{sec:decomp}
Under the assumption of Poisson arrival processes and exponential service-time distributions, our OQN is a
Markov model, called a Jackson network, which is easy to analyze, primarily because the steady-state distribution of the queue lengths has a product form; i.e.,
the steady-state queue lengths are independent geometric random variables,
just as if each queue were independent $M/M/1$
queues.
Motivated by that product-form property of Markov OQNs, decomposition approximations for non-Markov OQNs have been widely investigated.
In this approach, the network is decomposed into individual single-server queues, and the steady-state queue length processes are assumed to be approximately independent.
For example, in \cite{K79} and \cite{WW83qna}
each queue is approximated by a $GI/GI/1$ model, where the arrival (service) process is approximated by a renewal process partially characterized by the mean and \textit{squared coefficient of variation} of an interarrival (service) time.

While the decomposition approximations do often perform well, it was recognized that dependence in the arrival processes of the internal flows can be a significant problem.  The approximation for superposition processes used in the QNA algorithm \cite{WW83qna} attempts to address the dependence.  Nevertheless, significant problems remained, as was dramatically illustrated by comparisons of QNA to model simulations in \cite{SW86,FSW89,SW90}, as discussed in \cite{WW95}.

To address the dependence in arrival processes, decomposition methods based on Markov Arrival Processes (MAPs) have been developed. The
MAP was introduced by Neuts \cite{N79}; see Ch. XI of \cite{A03}.  Since it is not a renewal process, so that it can model the autocorrelation in the arrival and service processes.
 Horv\'{a}th et al. \cite{HHT10} approximated each station by a $\hbox{MAP}/\hbox{MAP}/1$ model.  Kim \cite{K11b,K11a} approximated each queue by a $\hbox{MMPP}(2)/GI/1$ model, where the arrival process is a Markov-modulated Poisson process with two states (a special MAP).

\subsubsection{Heavy-Traffic Limit Approximations}\label{sec:HT}

The early decomposition approximation in \cite{WW83qna} drew heavily on the central limit theorem (CLT) and heavy-traffic (HT) limit theorems.
Approximations for a single queue follow from \cite{IW70a,IW70b}.
With these tools, approximations for general point processes and arrival processes were developed in \cite{WW82point,WW84dep}.
Heavy-traffic approximation of queues with superposition arrival processes in \cite{WW85} helped capture the impact of dependence in such queues.

Another approach is to apply heavy-traffic limit theorems for the entire network.
Such HT limits were established for feedforward OQNs in \cite{IW70a,IW70b} and Harrison \cite{H73,H78}, and then for general OQNs by Reiman \cite{R84}.
These works showed that the queue length process converges to a multidimensional reflected Brownian motion (RBM) as every service station approaches full saturation simultaneously.

These general heavy-traffic results for OQNs lead to approximations 
using the
limiting RBM processes.
 The QNET algorithm in Harrison and Nguyen \cite{HN90} provides such an approximation.
Theoretical and numerical analysis of the stationary distribution of the multi-dimensional RBM was studied in \cite{HR81,HW87b,DH92}.

As a crucial step of the QNET algorithm, Dai and Harrison \cite{DH92} proposed a numerical algorithm to calculate the steady-state density of an RBM, but it require considerable computation time.
The computational accuracy of that algorithm improves as the number of iteration $n$ grows, and the author's there note that $n = 5$ generally gives satisfactory answers.
For a OQN with $d$ stations, the computational complexity is $O(d^{2n})$, see Section 6 of \cite{DH92}.
A further limitation is that the underlying theorem is for a sequence of OQNs in which the associated sequences
of traffic intensities at all queues approaches the critical value.
For practical application to large-scale systems or small systems with a wide range of traffic intensities, hybrid methods that combine a decomposition approximation and heavy-traffic theory were proposed in Reiman \cite{R90} and Dai et al. \cite{DNR94}.  The version \cite{DNR94} has been shown to be remarkably effective, but it requires the numerical solution
of RBMs.

\subsubsection{Robust Queueing Approximations}\label{sec:RQ}
Recently, a novel Robust Queueing (RQ) approach to analyze queueing performance in single-server queues has been proposed by Bandi et al. \cite{BBY15}.
The key idea in RQ is to replace the underlying probability law by a suitable uncertainty set, and analyze the
(deterministic) worst case performance.
The authors there relied on the discrete-time Lindley's recursion to characterize the customer waiting times as a supremum over partial sums of the interarrival times and service times.
Uncertainty sets for the sequence of partial sums are proposed based on central limit theorem and two-moment partial traffic descriptions of the arrival process and service process.

Although the general RQ idea is simple and good, there remain challenges in identifying proper uncertainty sets and making connection to the original queueing system.  These challenges were addressed
in \cite{WY18RQ1, WY19TVRQ}, which lay the foundation for this paper.
In \cite{WY18RQ1}, the authors proposed a new non-parametric RQ formulation for approximating the continuous-time workload process in a single-server queue, and proved asymptotic exactness of their approximations light and heavy traffic.
We briefly review this new RQ formulation in Section \ref{sec:RQ_review}.

\subsubsection{Non-Parametric Traffic Descriptions}\label{secNon-par}
As a trade-off for mathematical tractability, all approximation approaches so far rely on incomplete traffic descriptions.
For example, the approximation approaches reviewed above can be characterized as parametric approaches,
typically involving only means and variances of random variables,
The general stochastic system is then mapped into one of a parametric family of highly structured models.
Such approaches rely on a discrete set of parameters as traffic descriptions and a key step is to understand how these parameters evolve in the network.

Another stream of research models the temporal dependence in the stochastic processes by non-parametric traffic descriptions.
In Jagerman et al. \cite{JBAM04}, the authors approximate a general stationary arrival process by a Peakedness Matched Renewal Stream (PMRS).
The key ingredient is the peakedness function, which is determined by the arrival point process and the first two moments of the service-time distribution; see \cite{LiW14} for additional discussion..
However, \cite{JBAM04} relied on a two-parameter approximation for the peakedness function of a stationary point process, where the parameters are estimated by simulation.
Similar non-parametric traffic descriptions has been studied in \cite{LH92,LH93,JBAM04}, but they only focus on single-station single-server queues.

We adopt a non-parametric approach to describe the arrival and service processes in an OQN.
Let $A$ be a \textit{stationary} counting process, e.g. the arrival counting process at a queue.
We partially characterize $A$ by its rate and its \textit{Index of Dispersion for Counts} (IDC), a function of non-negative real numbers $I_A: \mathbb{R}^+ \rightarrow \mathbb{R}^+$ defined as in \S 4.5 of  \cite{CL66},
\begin{equation}\label{def: IDC}
I_A (t) \equiv \frac{\mathrm{Var}(A(t))}{E[A(t)]}, \quad t \ge 0.
\end{equation}
A reference case is the Poisson process, where the IDC is a constant function $I_{A}(t) \equiv 1$.
As regularity conditions, we assume that $E[A(t)]$ and $\mathrm{Var}(A(t))$ are finite for all $t\ge 0$.
For renewal processes, it suffices to assume that the inter-renewal time distribution have finite second moment.

Being a function of time $t$, the IDC captures the variability in a point process over all timescales.
The IDC encodes much more information about the underlying process
than traditional parametric descriptions.
The RQ algorithm in \cite{WY18RQ1} established a bridge between the IDC traffic description and the performance measures in a single-server queue.

With the aid of the HT limits established in \cite{WY18Dep,WY20flow},
we now develop a network calculus to characterize the IDCs of the customer flows in an OQN.
Similar non-parametric traffic descriptions have been studied in \cite{LH92,LH93,JBAM04}, but they focused on single-server queues.
To the best of our knowledge, we are the first to study the non-parametric traffic descriptions in a network setting.


\subsubsection{The Overall Robust Queueing Network Analyzer (RQNA)}\label{sec:RQNA}

We exploit the remarkably strong connection between the arrival IDC and the normalized workload in a single-server queue.
This connection, was first exposed by
Fendick and Whitt \cite{FW89}, but they did not produce
the systematic approximations we obtained through robust queueing in \cite{WY18RQ1}.
We advance that approach further by showing that all these approximations can be combined to produce a \textit{Robust Queueing Network Analyzer} (RQNA).

Our method is a decomposition approximation, because
the algorithm decomposes the network into individual $G/G/1$ models, where the arrival process and service process at each queue is partially specified by its rate and IDC, defined in (\ref{def: IDC}).
As in other decomposition methods, three network operations become essential:
first, the \textit{departure operation} as customers flow through a service station and an arrival process turns into a departure process;
second, the \textit{splitting operation} as a departure process split into multiple sub-processes and feed into different subsequent queues; and
third, the \textit{superposition operation} as departure flows from different queues combine together and feed into a queue.

In Section \ref{sec: IDC framework}, we discuss a set of linear equations, which we refer to as the \textit{IDC equations}, to describe the combined effect of these three network operations.
These IDC equations are derived from the HT limits in \cite{WY18Dep,WY20flow}.
We discuss the remaining technical details in the online supplement.
The IDC of the total arrival flows at each queue is approximated by the solution to the IDC equations.
The RQ algorithm (\ref{eqn: RQ approx}) is then applied to generate approximations of the mean steady-state performance measures at each $G/GI/1$ queue in the network.

The RQNA algorithm has a remarkably concise analytical formulation,
given in (\ref{eqn: RQ approx}) and (\ref{eqn: IDC equations}) below, which makes it easy to implement.
If we apply the feedback elimination in Section \ref{sec: feedback elimination},
 then the computational complexity is $O(K)$ or $O(K^2)$, where $K$ is the number of stations in the system.  We also conduct simulation experiments to evaluate the effectiveness of the new RQNA and compare it to previous algorithms in \cite{WW83qna,HN90,DNR94,HHT10}.
Our experiments indicate that RQNA performs as well or better than previous algorithms.

\subsection{Organization}
The rest of the paper is organized as follows.
In \S \ref{sec:review_IDC_RQ} we define the indices of dispersion, discuss the connection between the index of dispersion for work and the mean steady-state workload, and briefly review the robust queueing algorithm for a single $G/GI/1$ queue.
In \S \ref{sec: IDC framework} we develop a framework for approximating the IDC's of the flows.
In \S \ref{sec: RQNA tree} we develop a relatively elementary version of the RQNA algorithm for tree-structured networks.
In \S \ref{sec: feedback elimination} we discuss feedback elimination.
In \S \ref{sec: RQNA alg} we present the full RQNA algorithm.
 In \S \ref{sec: numerical}  we discuss numerical experiments.
We present additional supporting material in the appendix, including more experimental results.


\section{The Indices of Dispersion and Robust Queueing}\label{sec:review_IDC_RQ}

In this section we provide brief reviews of the IDC function in (\ref{def: IDC}) and the robust queueing algorithm
from \cite{WY18RQ1}.
In \S \ref{sec:def_index} we define another continuous-time index of dispersion: the Index of Dispersion for Work (IDW).
We discuss a useful decomposition of the IDW and its connection to the IDC and the mean steady-state workload.
In \S \ref{sec:RQ_review} we review the RQ algorithm from \cite{WY18RQ1}, which links the IDW to approximations of the steady-state queueing performance.

\subsection{The Indices of Dispersion}\label{sec:def_index}

Consider a general single-server queue with a general arrival process $A$, i.e. $A(t)$ counts the total number of arrival in the time interval $[0,t]$.  We assume that $A$ is a stationary point process; see \cite{DV08II,S95}.
The IDC defined in (\ref{def: IDC}) is a continuous-time function associated with $A$.
Being the variance function scaled by the mean function, the IDC exposes the variability over time, independent of the scale.
For this reason, the IDC can be viewed as a continuous-time generalization
of the squared coefficient of variation ($scv$, variance divided by the square of the mean) of a nonnegative random variable. The IDC captures the way that the covariance in a point process changes over time, which extends the natural practice of including lag-$k$ covariances in modeling the dependence in a point process.

The reference case is a Poisson arrival process, for which $I_a (t) = 1$, $t \ge 0$.
However, for general arrival processes, the IDC is more complicated.
Even the IDC for a determinsitic $D$ arrival process is complicated, because the IDC is for the stationary version of the arrival process, which lets the initial point be uniformly distributed over the constant interarrival time.
Much of this paper is devoted to the analysis and approximation of the IDC for the arrival process at each station of the OQN.

\begin{remark}\label{rmRenewalIDC}
\em{
(Advantage of using the IDC)
In \cite{WY19Index}, we discussed the advantage of using the IDC in queueing approximations.
Theorem 2.1 there showed that, for renewal processes, the inter-renewal time distribution can be fully recovered from the rate and the IDC of the associated equilibrium renewal process.
Thus, when considering an OQN with renewal external arrival processes and i.i.d. service times, the network is fully characterized by the rate and IDC's of the external arrival processes and the service processes (both regarded as stationary processes).
Hence, the IDC function encodes much more information about the underlying distribution than traditional traffic descriptions. \hfill $\square$
}
\end{remark}

\begin{remark}\label{rmScalingConvention}
\em{
(Time scaling convention)
In \cite{WY18RQ1} we defined the IDC and IDW in terms of rate-$1$ processes, so that the actual rate of the process had to be inserted as part of the time argument.  In contrast, here as in \cite{WY18Dep} we let the underlying processes $A$ and $Y$ have any given rate, so no further scaling is needed.  That changes the formulas for the IDC of a superposition process, e.g., compare (36) of \cite{WY18RQ1} to \eqn{superposition} here.   To illustrate the idea, consider $A(t)$ with rate-$1$ and $A_{\lambda}(t) \equiv A(\lambda t)$ with rate-$\lambda$.  Let $I_{A}(t)$ denote the IDC of $A(t)$, then we have $I_{A_{\lambda}}(t) \equiv \Var(A(\lambda t))/E[A(\lambda t)] = I_{A}(\lambda t). $ \hfill $\square$
}
\end{remark}

Now, consider a general sequence of service times $\{V_i: i\ge 1\}$, where $V_i$ is the service requirement of the $i$-th customer.
Let
\begin{equation}\label{def: cumulative work}
Y(t) \equiv \sum_{i = 1}^{A(t)} V_i
\end{equation}
denote the \textit{cumulative work input process}.
This process connects to the workload of a single-server queue by (\ref{c2}) and (\ref{c3}) below.

Paralleling the IDC, the \textit{Index of Dispersion for Work} (IDW) describes the variability associated with the cumulative input process $Y$ in (\ref{def: cumulative work}).
The IDW is defined as in (1) of \cite{FW89} by
\begin{equation}\label{def: IDW}
I_w (t) \equiv \frac{\Var(Y (t))}{E[V_1] E[Y (t)]}, \quad t \ge 0.
\end{equation}
The IDW captures the cumulative variability of the total service requirement brought to the system as a function of time $t$, which is a key component of the new RQ approximation in \cite{WY18RQ1} as we review in \S \ref{sec:RQ_review}.

Since we are interested in the steady-state performance of the OQN, we assume that the processes $A$ and $Y$ have stationary increments.
Given that arrival process and service times have constant determined rates, the mean functions $E[A(t)]$ and $E[Y(t)]$ are linear in time.
Hence, much of the remaining behavior of the $A$ and $Y$ is determined by the variance-time function or index of dispersion.
We are interested in the variance-time {\em function}, because it captures the dependence through the covariances;
the processes $(A,Y)$ have independent increments for the $M/GI/1$ model, but otherwise not.

To connect the IDC to the IDW, consider the special case where the service times $V_i$ are i.i.d, independent of the arrival process $A(t)$.
The conditional variance formula gives a useful decomposition of the IDW
\begin{equation}\label{eqn: IDW decomposition}
I_w(t) = I_a(t) + c_{s}^2, \quad t \ge 0,
\end{equation}
where $c_s^2 = \Var(V_i)/E[V_i]^2$ is the scv of the service-time distribution.

\subsubsection{The IDW and the Mean Steady-State Workload}\label{secLink}

The IDC and IDW are important because of their close connection to the mean steady-state workload $E[Z_{\rho}]$.
Here we make the performance measure explicitly depend on the traffic intensity $\rho$ to expose the joint impact of dependence in flows and the traffic intensity on it.
Under regularity conditions, the workload $Z(t)$ converges to the steady-state workload $Z_{\rho}$ as $t$ increases to infinity.
In \cite{FW89} it was shown that the IDW $I_w$ is intimately related to a scaled mean workload $c^2_Z (\rho)$, defined
by
\beql{z4}
c^2_{Z} (\rho) \equiv \frac{E[Z_{\rho}]}{E[Z_{\rho}; M/D/1]},
\eeq
where $E[Z_{\rho}; M/D/1]$ is the mean steady-state workload in a M/D/1 model given by
\begin{equation}\label{workload MD1}
E[Z_{\rho}; M/D/1] = \frac{E[V_1] \rho}{2(1-\rho)}.
\end{equation}
As (\ref{workload MD1}) suggests, the mean steady-state workload converges to $0$ as $\rho \downarrow 0$ and diverges to infinity as $\rho \uparrow 1$.
The normalization in \eqn{z4} exposes the impact of variability separately from the traffic intensity.

In great generality as discussed in \cite{FW89}, we have
\begin{equation}\label{gen}
c^2_{Z} (0) = 1 + c_s^2 = I_{w} (0) \qandq c^2_{Z} (1) = c_{A}^2 + c_s^2 = I_{w} (\infty),
\end{equation}
where $c_{A}^2$ is the asymptotic variability parameter, i.e., the normalization constant in the
functional central limit theorem (FCLT) for the arrival process.  For a renewal process, $c_A^2$ coincides with the scv $c_a^2$ of an interarrival time.
The reference case is the classical $M/GI/1$ queue, for which we have
\[c^2_{Z} (\rho) = 1 + c_s^2 = I_{w} (t) \qforallq \rho, \, t, \quad 0 < \rho < 1, \, t \ge 0.\]

The limits in (\ref{gen}) imply that, when $c_{A}^2$ is not nearly $1$, $c_Z^2 (\rho)$ varies significantly as a function of $\rho$.  Hence, the impact of the variability in the arrival process upon the queue performance clearly depends on
the traffic intensity.  This important insight from \cite{FW89} is the starting point for our analysis.
In well-behaved models, $c_Z^2 (\rho)$ as a function of $\rho$ and $I_w (t)$ as a function of $t$
tend to change smoothly and monotonically between those extremes, but OQNs can produce more complex behavior when both the traffic intensities at the queues and the levels of variability in the arrival and service processes at different queues vary; e.g., see the examples for queues in series in \S\S 5.2, EC.8.2 and EC8.3 of \cite{WY18RQ1}.

\subsection{Robust Queueing for Single-Server Queues}\label{sec:RQ_review}
In this section, we review the RQ algorithm for single-server queues and discuss approximations for other performance measures obtained as a result.
The RQ algorithm serves as a bridge between the IDC of the arrival process and the approximations of the performance measures.
In particular, as in (\ref{eqn: RQ approx}), the RQ algorithm generates approximation of the steady-state workload for any queue using the IDC of the total arrival process at that queue.  

Consider the $G/GI/1$ queue, where the arrival process is a stationary and ergodic point process and the service times are i.i.d., independent of the arrival process.
We assume that the arrival process $A$ is partially characterized by the arrival rate $\lambda$ and the IDC $I_a$ defined in (\ref{def: IDC}).
For a stationary point process, we always have $E[A(t)] = \lambda t$; see \S 2.7 of \cite{S95}.
We further assume that the service time distribution has finite mean $1/\mu$ (and thus rate $\mu$) and scv $c_s^2$.
 We also assume that $\rho \equiv \lambda/\mu < 1$ for model stability.
Let $Z$ be the steady-state workload in the $G/GI/1$ model.
The RQ algorithm provides approximation for $E[Z]$ with $(\lambda, I_a,\mu, c_s^2)$ as input data.

To obtain the RQ algorithm, we start with a reverse-time construction of the workload process as in \S 3 of \cite{WY18RQ1}.
Define the net-input process $N(t)$ as
\beql{c2}
N (t) \equiv Y (t) - t, \quad t \ge 0.
\eeq
Then the workload at time $t$, starting empty at time  $0$, is obtained from
the reflection map $\Psi$ applied to $N$, i.e.,
\beql{c3}
Z = \Psi(N) (t) \equiv  N (t) - \inf_{0 \le s \le t}{\{N (s)\}}, \quad t \ge 0.
\eeq
With a slight abuse of notation,
let $Z(t)$ be the workload at time $0$ of a system that started empty at time $-t$.  Then $Z(t)$ can be represented as
\beql{z1}
Z(t) \equiv \sup_{0 \le s \le t}{\{N (s)\}}, \quad t \ge 0,
\eeq
where $N$ is defined in terms of $Y$ as before, but $Y$ is interpreted as the total work in service time to enter over the interval $[-s,0]$.
That is achieved by letting
$V_k$ be the $k^{\rm th}$ service time indexed going backwards from time $0$ and
$A(s)$ counting the number of arrivals in the interval $[-s,0]$.

The workload process $Z(t)$ defined in \eqn{z1} is nondecreasing in $t$ and hence necessarily converges to
a limit $Z$.
For the stable stationary $G/GI/1$ model, $Z$ corresponds to the steady-state workload and satisfies $P(Z < \infty) = 1$; see \S 6.3 of \cite{S95}.

In the ordinary stochastic queueing model, $N(s)$ is a stochastic process and hence $Z(t)$ is a random variable.
However, in Robust Queueing practice, $N(s)$ is viewed as a deterministic instance drawn from a pre-determined uncertainty set $\sU$ of input functions, while the workload $Z^*$ for a Robust Queue is regarded as the worst case workload over the uncertainty set, i.e.
\[Z^* \equiv \sup_{\tilde{N} \in \sU} \sup_{x \ge 0} \{\tilde{N}(x)\}.\]
Following the setting from \citep{WY18RQ1}, we adopt the following uncertainty set motivated from central limit theorem (CLT)
\begin{align}
\sU  & \equiv  \left\{\tilde{N}: \RR^{+} \to \RR: \tilde{N} (s)\le E[N (s)] +  b\sqrt{\mathrm{Var}(N (s))}, \, s \ge 0 \right\},
\end{align}
where $N(t)$ is the net input process associated with the stochastic queue, so
\begin{align*}
  E[N(t)] &  = E[Y(t) - t] = \rho t - t,\\
  \mathrm{Var}(N(t)) & = \mathrm{Var}(Y(t)) = I_w(t)E[V_1]E[Y(t)] =  (I_a(t) + c_s^2)\rho t/\mu.
\end{align*}
The RQ approximation based on this partial model characterization is
\begin{align}
E[Z_{\rho}] \approx Z_{\rho}^{*} \equiv \sup_{\tilde{N}_{\rho} \in \sU_{\rho}} \sup_{x \ge 0} \{\tilde{N}(x)\}
 & = \sup_{x \ge 0}{\{ -(1-\rho)x + b\sqrt{\rho x (I_a (x) + c_s^2)/\mu}\}}, \label{eqn: RQ approx}
\end{align}
which follows Theorem 2 of \cite{WY18RQ1} and (\ref{eqn: IDW decomposition}).
Notice that the approximation in \eqn{eqn: RQ approx} is directly a supremum of a real-valued function,
and so can be computed quite easily for any given 4-tuple $(\lambda, I_a, \mu, c_s^2)$.

Theorem 5 in \cite{WY18RQ1} states that the RQ algorithm gives asymptotically exact values of the mean steady-state workload in both light-traffic and heavy-traffic limits.  Through extensive simulation experiments, it has been found that the mean steady-state workload $E[Z]$ can be well approximated by the IDW-based RQ algorithm.


\begin{remark}\label{Rmk: finite sets}
\em{
(Finitely many times)
Even though the IDC $I_{a} (t)$ and the supremum over $x$ in \eqn{eqn: RQ approx} are defined
over the interval $[0, \infty)$,
we work with a finite subset of times.  In \S 3 of the appendix, we develop an algorithm for computing $I_{a} (t)$
for any given $t$.
In particular, we use the supremum over a grid equally spaced in logarithm time scale in \eqn{eqn: RQ approx}.
To justify the logarithm scale, we remark that the IDC of common point processes converges exponentially fast to a constant, as the time $t$ increases.  In particular, this holds for Markov arrival processes, which includes hyper-exponential renewal process, Erlang renewal process, and Markov modulated Poission Process as special cases;
e.g., see Ch. XI of \cite{A03}.
\hfill $\square$
}
\end{remark}

\begin{remark}\label{Rmk: stationary}
\em{
(Continuous-time stationarity)
We emphasize that, in the RQ formulation, it is essential to use the continuous-time stationary version of the IDC in \eqn{def: IDC} and the IDW in \eqn{def: IDW}, instead of their discrete-time Palm stationary versions; see \cite{S95} for a comprehensive discussion.    The continuous-time stationary IDC we use here yields asymptotically correct light-traffic limit, whereas the Palm stationary IDC does not; see \S 5.2 of \cite{WY18RQ1}.  \hfill $\square$
}
\end{remark}

\begin{remark}\label{Rmk:other_perf}
\em{
(Queue length and waiting time)
Approximations for other steady-state performance measures can be obtained by applying exact relations for
the $G/GI/1$ queue that follow from Little's law $L = \lambda W$ and its generalization $H = \lambda G$;
e.g., see \cite{WW91rev} and
Chapter X of \cite{A03} for the $GI/GI/1$ special case.
Let $W, Q$ and $X$ be the steady-state waiting time,
queue length and the number in system (including the one in service, if any).  By Little's law,
\begin{align*}\label{LL}
E[Q] & = \lambda E[W] = \rho E[W] \qandq \nonumber \\
E[X] & =  E[Q] + \rho = \rho (E[W] + 1).
\end{align*}
By Brumelle's formula \cite{B71} or $H = \lambda G$, (6.20) of \cite{WW91rev},
\[E[Z] = \rho E[W] + \rho\frac{E[V^2]}{2\mu} = \rho E[W] + \rho\frac{(c_s^2 + 1)}{2\mu}.\]
Hence, given an approximation $Z^*$ for $E[Z]$, we can use the approximations
\begin{align*}
E[W] & \approx \max\{0, Z^*/\rho - (c_s^2 + 1)/2\mu \} \qandq \nonumber \\
E[Q] & \approx \lambda E[W].
\end{align*}
}
\end{remark}

\begin{remark}\label{rmNet}
\em{
(Network performance measures)
So far we only have discussed the performance measures for a single station.  The total network performance measures, on the other hand, can also be derived.  For example, the expected value of the total sojourn time  $T_{i}^{\mathrm{tot}}$, i.e. the time needed to flow through the queueing network for a customer that enters the system from station $i$, is easily estimated from the obtained mean waiting time at each station.  Assuming Markov routing with routing matrix $P$, a standard argument from discrete time Markov chain theory gives the mean total number of visits  $\xi_{i,j}$ to station $j$ by a customer entering the system at station $i$ as
\[\xi_{i,j} = \left((I - P)^{-1}\right)_{i,j},\]
where $(I - P)^{-1}$ is the fundamental matrix of a absorbing Markov chain.
Hence, the mean steady-state total sojourn time $E[T^{\mathrm{tot}}_{i}]$ is approximated by
\begin{equation}\label{eqn: total sojourn}
  E[T^{\mathrm{tot}}_{i}] \approx \sum_{j = 1}^K \xi_{i,j}(E[W_j] + 1/\mu_j).
\end{equation}
In real world applications, customers often experiences non-Markovian routing, where routes are customer-dependent.
For ways to represent those scenarios and convert them (approximately) to the current framework, see \S 2.3 and \S 6 of \cite{WW83qna}. \hfill $\square$
}
\end{remark}

\section{Approximating the IDCs of the Network Flows}\label{sec: IDC framework}

In the i.i.d. service time setting, the IDW reduces to the arrival IDC plus the service scv as in (\ref{eqn: IDW decomposition}).
To generalize the RQ algorithm in \S \ref{sec:RQ_review} into a RQNA algorithm for networks,  the main challenge is developing a successful approximation for the IDC of the total arrival flow at each queue.

In this section we develop a framework for approximating the IDCs of the network flows in the OQN, including the total arrival flows.
We start in \S \ref{secData} by reviewing the OQN model and the required model data for the RQNA algorithm.
We review the standard traffic rate equations in \S \ref{secRateEqns} and develop the new IDC equations in \S \ref{secVarEqns}.

\subsection{The OQN Model}\label{secData}
\subsubsection{The Model Primitives} \label{sec: model primitives}
We consider a network of $K$ queues. Each queue has a single server, unlimited waiting space and provides service in order of arrival.

For each queue $i$, $1 \le i \le K$, we have an external arrival process $A_{0,i} \equiv \{A_{0,i} (t): t \ge 0\}$.  Each external arrival process $A_{0,i}$ is assumed to be a simple (no batches) stationary and ergodic point process with finite rate $\lambda_{0,i}$ and finite second-moment process $E[A_{0,i}^2 (t)]$.
We assume that all these external arrival processes, as well as the service and routing processes, are mutually independent.

For each individual queue, we assume that the service times are i.i.d.  Let $V_i^{l}$ denote the service requirement of the $l$-th customer at queue $i$, which we assume to be distributed according to cdf $G_i$ with finite mean $1/\mu_i$ and scv $c_{s,i}^2$.  Let the associated service renewal counting process be $S_i \equiv \{S_i (t): t \ge 0\}$, where
\begin{equation}\label{eqn: service renewal}
S_i(t) = \max\left\{n\leq 0: \sum_{l = 1}^n V_i^l \leq t\right\}, \quad t \ge 0.
\end{equation}

 We assume that departures are routed from node to node and out of the network by Markovian routing, which is independent of the arrival and service processes.  We assume that each arrival eventually leaves w.p.1. Let $p_{i,j}$ denote the probability that a departure from node $i$ is routed to node $j$. Let $P \equiv \{p_{i,j}: 1\le i,j \le K\}$ be the (substochastic) routing matrix.  Furthermore, let $p_{i,0} \equiv 1- \sum_{j} p_{i,j}$ denote the probability that a customer departs the system after completing service at from node $i$.

\subsubsection{The IDC's of the Flows}\label{sec: IDC flows}
In order to apply the RQ algorithm, our primary focus here is to analyze and approximate the IDC's of the customer flows in a OQN.  
The flows can be separated into two groups, the \textit{external flows} and the \textit{internal flows}.  The external flows are the flows associated with the model primitives in \S\ref{sec: model primitives}.  For external arrival process $A_{0,i}$, we let $I_{a,0,i} \equiv \{I_{a,0,i}(t): 0 \le t \le \infty\}$ denote the its IDC, as defined in \eqn{def: IDC}.
For service flows, let $I_{s,i} \equiv \{I_{s,i}(t); 0 \le t \le \infty\}$ be the IDC of the stationary renewal process associated with (\ref{eqn: service renewal}).  For the case of renewal process, we necessarily have $I_{s,i} (\infty) = c_{s,i}^2$. We assume that the IDC's $I_{a,0,i}$ and $I_{s,i}$ are continuous functions with limits at $0$ and $+\infty$.

The IDC's of the external flows forms an important part of the model input of our RQNA algorithm.  In particular, we assume that we are given $(\lambda_{0,i}, I_{a,0,i}, \mu_i, I_{s,i})$ for each queue $i$ and the routing matrix $P$.

In practice, the IDC of the external flows can be specified by one of the following ways.  First, for renewal processes, it suffices to specify the interrenwal-time cdf; then
the associated IDC can be computed from the cdf as indicated
 in \S 3.1 of the appendix. 
 Second, if we are only given the first two moments, then we can fit a convenient cdf to these parameters as indicated in \S 3 of \cite{WW82point}, and use the corresponding IDC.  Third, if we are only give sample data of the process, then we apply the numerical algorithm in \S 3.2 of the appendix
 to estimate the rate and IDC of the process.

To implement our IDC approximations, we develop approximations for the IDC's of the internal flows.
We use the following notation:  Let $A_{i}$ denote the total arrival process at queue $i$ and let $I_{a,i}$ be the associated IDC; let $D_{i}$ denote the departure process at queue $i$ and let $I_{d,i}$ be the associated IDC;
and let $A_{i,j}$ denote the departing customer flow from queue $i$ that are routed to queue $j$ and let $I_{a,i,j}$ be the associated IDC.

\subsection{The Traffic Rate Equations and Traffic Intensities}\label{secRateEqns}

Let
$\lambda \equiv (\lambda_{1}, \dots, \lambda_{K})$ be the effective (total) arrival rate vector.  We use the same traffic rate equations as in a Jackson network to determine $\lambda$.
Then  $\lambda_{i,j} \equiv \lambda_{i} p_{i,j}$ is the rate of the internal arrival flow $A_{i,j}$.  Recall that $\lambda_0 \equiv (\lambda_{0,1}, \dots, \lambda_{0,K})$ is the external arrival rate vector, then the traffic-rate equations are
\begin{equation}\label{traffic_rate_equations}
\lambda_i = \lambda_{0,i} + \sum_{j = 1}^K \lambda_{j,i} = \lambda_{0,i} + \sum_{i = 1}^K \lambda_{j} p_{j,i}, \quad 1 \le i \le K,
\end{equation}
or in matrix form
\[(I - P')\lambda = \lambda_0,\]
where $I$ denotes the $K\times K$ identity matrix.
We assume that $I-P'$ is invertible; i.e., we assume that all customers eventually leave the system.  The condition for the invertibility of $I-P'$ to hold is well known, e.g. in Theorem 3.2.1 of \cite{KS76}.
Hence, the vector of internal arrival rates is given by
\beql{traffic_rate_equations_sol}
\lambda = (I - P')^{-1}\lambda_0.
\eeq

Then the traffic intensity at queue $i$ is defined as usual by $\rho_i \equiv \lambda_i/\mu_i$.
We assume that $\rho_i < 1$ for all $i$ to ensure that the OQN is stable.

\subsection{The Traffic Variability Equations}\label{secVarEqns}
In this section, we develop a set of IDC equations to solve for the approximations of the IDC's of the internal flows.  The IDC of the total arrival process at each queue is then converted into approximations of the performances measures as in \S \ref{sec:RQ_review}.

As in other decomposition methods, three network operations are essential: the departure operation (flow through a queue), the splitting operation (devide a flow into several sub-flows) and the superposition operation (combining multiple flows). We develop IDC equations that reveal (approximately) how the IDC's evolve under each network operation.


\subsubsection{The Departure Operation}\label{sec: review HTDep}

The IDC of the stationary departure process has been studied in \S 6.2 of \cite{WY18Dep}.  We briefly review the departure IDC equation, see \S 5.1 of the appendix 
for more details.

We approximate the IDC $I_{d,i}$ by a convex combination of the
arrival IDC $I_{a,i}$ and the service IDC $I_{s,i}$. In particular,
\beql{depIDC}
  I_{d,i}(t) \approx w_{i}(t) I_{a,i} (t) + (1-w_{i}(t)) I_{s,i} (\rho_i t), \quad t \ge 0.
\eeq
The weight function $w_{i}$ is defined as
 \beql{depIDC2}
  w_{i}(t) \equiv w^{*} \left((1-\rho_i)^{2} \lambda_i t/\rho_i c_{x,i}^2\right), \quad t \ge 0,
\eeq
where $c_{x,i}^2 \equiv c_{a,i}^2 + c_{s,i}^2$ and $c_{a,i}^2 = I_{a,i} (\infty)$ and the {\em canonical weight function} $w^*$ is
\begin{equation}\label{weight}
w^{*}(t) = \frac{1}{2t}\left(\left(t^2 + 2t -1\right)\left(1 - 2\Phi^c(\sqrt{t})\right) + 2\phi(\sqrt{t})\sqrt{t}\left(1+t\right) - t^2 \right)
\end{equation}

Note that there is a change of notation between (\ref{depIDC}) here and (74) in \cite{WY18Dep}.  In particular, we have $I_{s,i} (\rho_i t)$ here instead of $I_{s,i}(t)$.  In \cite{WY18Dep}, we worked with a single-server queue and assumed that $I_{s,i}(t)$ is the IDC associated with the rate-$\lambda_i$ service process.  However, when considering a OQN here, it is natural to work with service IDC that associated with the service rate $\mu_i$.   These two approaches are equivalent, as we observed in Remark \ref{rmScalingConvention}.
Given that the given stationary service process has rate $\mu_i$, we convert it to rate $\lambda_i$ by considering $I_{s,i} (\rho_i t)$.

\begin{remark}\label{rmDep}
{\em
(Parallel to QNA in \cite{WW83qna}.)
The convex combination in the approximation \eqn{depIDC} is reminiscent of the
convex combination for variability parameters in (38) of \cite{WW83qna}, i.e.,
\beql{depscv}
c_{d,i} \approx (1 - \rho_i^2)c_{a,i}^2 + \rho_i^2 c_{s,i}^2,
\eeq
which corresponds to a stationary-interval approximation, as discussed in \cite{WW82point,WW83qna,WW84dep}.

Similar behavior can be seen in approximation \eqn{depIDC}.  In particular,
the canonical weight function $w^*$ in \eqn{weight} is a monotonically increasing function with $w^*(0) = 0$ and $w^*(\infty) = 1$.  By the definition of $w_{i}(t)$,
we see that for each $t$, (\ref{depIDC}) places less weight on $I_{a,i} (t)$ and
more weight on $I_{s,i} (t)$
as $\rho_i$ increases.
This makes sense intuitively, because the
queue should be busy most of the time as $\rho_i$ increases toward $1$.  Thus departure times
tend to be minor variations of service times.  In contrast, if $\rho_i$ is very small, then the queue acts only as a minor perturbation of the arrival process.

However, \eqn{depIDC2} reveals a more subtle interaction between $\rho_i$ and the variability of the departure process over different time scales. \hfill $\square$
}
\end{remark}

\subsubsection{The Splitting Operation}\label{sec: splitting}

To treat splitting, we write the split process $A_{i,j}$ as a random sum. Let
 $\theta^{l}_{i,j} = 1$ if the $l$-th departure from queue $i$ is directed to queue $j$, and let $\theta^{l}_{i,j} = 0$ if otherwise.
Then observe that
\[A_{i,j} (t) = \sum_{l = 1}^{D_i (t)} \theta^{l}_{i,j}, \quad t \ge 0.\]
We apply the conditional-variance formula to write the variance $V_{a,i,j}(t) \equiv \mathrm{Var}(A_{i,j}(t))$ as
\begin{eqnarray}\label{condVar}
V_{a,i,j}(t) & = & E[\Var(A_{i,j} (t)|D_i (t))] + Var(E[A_{i,j} (t)|D_i (t)]).
 \end{eqnarray}

With the Markovian routing we have assumed, the routing decisions at each queue at each time are i.i.d. and
independent of the history of the network.
As a consequence, for feed-forward queueing networks, we can deduce that the collection of all routing decisions
made at queue $i$ up to time $t$ is independent of $D_i (t)$.  For the case in which independence holds, we can apply \eqn{condVar} to express $V_{a,i,j}(t)$
in terms of
the variance of the departure process, $V_{d,i}(t) \equiv \mathrm{Var}(D_{i}(t))$;
in particular,
\beql{splitVar}
V_{a,i,j}(t) = p^2_{i,j} V_{d,i}(t)  + p_{i,j}(1-p_{i,j})\lambda_i t,
\eeq
or, equivalently, since $E[D_{i}(t)] = \lambda_i t$ and $E[A_{i,j}(t)] = p_{i,j}\lambda_i t = p_{i,j}E[D_{i}(t)]$,
\begin{equation}\label{splitting}
I_{a,i,j}(t) = p_{i,j} I_{d,i}(t)  + (1-p_{i,j}).
\end{equation}
The formula \eqn{splitting} is an initial approximation, which parallels the
  approximation used for splitting in (40) of \cite{WW83qna}, i.e.,
  $c_{a,i,j}^2 = p_{i,j}c_{d,i}^2 + (1-p_{i,j})$.

However, the independence assumption will not hold in the presence of customer feedback, in which case there is a complicated dependence.
we develop a more general formula to improve the approximation in general OQNs.

For that purpose, we apply the FCLT for split processes in \S 9.5 of \cite{WW02} and the heavy-traffic limit theorems in \cite{WY20flow}.  We give the detailed derivation in \S 5.2 of the appendix. 

Based on that heavy-traffic analysis, we propose the splitting IDC equation as
\begin{equation}\label{eqn: splitting equation}
I_{a,i,j}(t) = p_{i,j}I_{d,i}(t) + (1-p_{i,j}) + \alpha_{i,j}(t).
\end{equation}
To account for the dependence, we include a correction term $\alpha_{i,j}$, defined as
\begin{eqnarray}\label{eqn: correction alpha}
\alpha_{i,j,\rho_i}(t) & \approx & 2\xi_{i,j} p_{i,j}(1-p_{i,j})w_{\rho_i}(t) \nonumber \\
& = & 2\xi_{i,j} p_{i,j}(1-p_{i,j})w^*((1-\rho_i)^{-2}\lambda_i t/(h(\rho_i) c_{x,i}^2)), \quad t \ge 0,
\end{eqnarray}
where $w_{\rho_i}(t)$ is the weight function for the departure IDC in \eqn{depIDC2},
$c_{x,i}^2$, $c_{a,i}^2$ and $c_{s,i}^2$ are also as in \eqn{depIDC2}, while
$\xi_{i,j}$ is the $(i,j)^{\rm th}$ entry of the matrix $(I-P')^{-1}$.

\subsubsection{The Superposition Operation}\label{sec: superposition}
In this section, we investigate the impact of the superposition operation on the IDC's.  To start, consider the case in which the individual streams are mutually independent.  In this case, we have
\[V_{a,i}(t)\equiv \mathrm{Var}(A_{i}(t)) = \mathrm{Var}\left(\sum_{j = 0}^{K} A_{j,i}(t)\right) = \sum_{j = 0}^{K} \mathrm{Var}(A_{j,i}(t)) ,\]
so that
\begin{equation}\label{superposition}
I_{a,i} (t) = \sum_{j = 0}^{K} (\lambda_{j,i}/\lambda_i) I_{a,j,i}(t),
\end{equation}
where $I_{a,j,i}(t) \equiv \mathrm{Var}(A_{j,i}(t))/E[A_{j,i}(t)]$.
Recall that \eqn{superposition} differs from (36) of \cite{WY18RQ1} because we are not assuming rate-$1$
processes in our definitons of the IDC; see Remark \ref{rmScalingConvention}.

While (\ref{superposition}) is exact when the streams are independent, it is not exact in general cases.  Even for feed-forward networks, we may have a stream that splits and then recombines later, which introduces dependence.

For dependent streams, the variance of the superposition total arrival process at queue $i$ can be written as
\[V_{a,i}(t) \equiv \textrm{Var}\left(\sum_{j = 0}^K A_{j,i}(t)\right) = \sum_{j = 0}^K \textrm{Var}\left(A_{j,i}(t)\right) + \beta_i(t)E[A_i(t)]\]
where $A_{0,i}$ denotes the external arrival process at station $i$,
\begin{equation}\label{eqn: superposition correction term - general}
\beta_i(t) \equiv \sum_{j \neq k} \beta_{j,i;k,i}(t), \qandq \beta_{j,i;k,i}(t)
\equiv \frac{\textrm{cov}\left(A_{j,i}(t), A_{k,i}(t)\right)}{E[A_i(t)]}.
\end{equation}
In terms of the IDC's, we have
\begin{equation}\label{eqn: superposition equation}
  I_{a_i}(t) = \sum_{j = 0}^{K} (\lambda_{j,i}/\lambda_i) I_{a_{j,i}}(t) + \beta_i(t).
\end{equation}

In general, an exact characterization of the correction term $\beta_{i}(t)$ is not available.  Thus, we again apply heavy-traffic limits in \cite{WY20flow}
to generate an approximation. Detailed derivation appears in \S 5.3 of the appendix. 

Assume without loss of generality that $\rho_j \ge \rho_i$.
From the heavy-traffic analysis,
we obtain the approximation
\begin{equation}\label{eqn: correction beta}
\beta_{j,i;k,i}(t) = \beta_{k,i;j,i}(t) \approx  (\zeta_{j,i;k,i}/\lambda_i) w^*((1-\rho_j)^{2}p_{j,i} \lambda_j t/\rho_j c_{x,j,i}^2),
\end{equation}
where $w^*$ is the weight function in \eqn{weight},
$c_{x,j,i}^2 = p_{j,i} c_{a,j}^2 + (1-p_{j,i}) + p_{j,i} c^2_{s,j}$ and $c_{a,j}^2$ is solved from the variability equations for the asymptotic variability parameters in (\ref{eqn: limiting variability equations}).  The constant $\zeta_{j,i;k,i}$ is defined as
\begin{equation}\label{eqn: zeta}
 \zeta_{j,i;k,i} = \nu_j' \left(\mathrm{diag}(c^2_{a,0,i} \lambda_i) + \sum_{l = 1}^K \Sigma_l\right) \nu_k + \nu'_k \Sigma_j e_i + \nu'_j\Sigma_k e_i,
\end{equation}
where $\nu_l \equiv p_{l,i}e'_l (I-P')^{-1}$ for $l = j,k$, $e_i$ is the $i$-th unit vector, $\mathrm{diag}(c^2_{a,0,i} \lambda_i)$ is the diagonal matrix with $c^2_{a,0,i} \lambda_i$ as the $i$-th diagonal entry,  $\Sigma_l$ is the covaraince matrix of the splitting decision process at station $l$ defined as $\Sigma_l \equiv (\sigma^l_{i,j})$ with $\sigma^l_{i,i} = p_{l,i}(1-p_{l,i})\lambda_l$ and $\sigma^l_{i,j} = -p_{l,i} p_{l,j}\lambda_l$ for $i\neq j$.


\subsection{The IDC Equation System}\label{sec: IDC equations}

We now assemble the building blocks into a system of linear equations (for each $t$) that describes the IDC's in the OQN.
Combining (\ref{depIDC}), (\ref{eqn: splitting equation}) and (\ref{eqn: superposition equation}), we obtain \textit{the IDC equations}.  These are equations that should be satisfied by the unknown IDCs. For $1 \le i \le K$, the equations are
\begin{align}
 I_{a,i}(t) & =  \sum_{j = 1}^{K} (\lambda_{j,i}/\lambda_i) I_{a,j,i}(t) + (\lambda_{0,i}/\lambda_i) I_{a,0,i}(t) + \beta_i(t),   \nonumber \\
 I_{a,i,j}(t)& =  p_{i,j}I_{d,i}(t) + (1-p_{i,j}) + \alpha_{i,j}(t),  \nonumber\\
 I_{d,i}(t) & =  w_i(t) I_{a,i} (t) + (1-w_i(t)) I_{s,i}(\rho_i t).\label{eqn: IDC equations - nonmatrix}
 \end{align}
The parameters $p_{i,j}$, $\lambda_{i,j}$ and $\lambda_i$ are determined by the model primitives in \S \ref{sec: model primitives} and the traffic rate equations in \S \ref{secRateEqns}.
The IDC's of the external flows $I_{a_{0,i}}(t)$ and $I_{s_i}(t)$ are assumed to be calculated via exact or numerical inversion of Laplace Transforms, or estimated from data.
The weight functions $w_i(t)$ is defined in (\ref{depIDC2}), which involves a limiting variability parameter $c_{x,i}^2 \equiv I_{a,i}(\infty) + c_{s,i}^2$.

To solve for the limiting variability parameters $I_{a,i}(\infty)$, we let $t\to \infty$ in (\ref{eqn: IDC equations - nonmatrix}) and denote $c^2_{a,i} \equiv I_{a,i}(\infty), c^2_{a,i,j} \equiv I_{a,i,j}(\infty)$ and $c^2_{d,i} \equiv I_{d,i}(\infty)$.   Furthermore, we define
\begin{align*}
c^2_{\alpha_{i,j}} &\equiv \alpha_{i,j}(\infty) = 2 \xi_{i,j} p_{i,j}(1 - p_{i,j}), \\
c^2_{\beta_i}  &\equiv \beta_i(\infty)  = \frac{2}{\lambda_i}\sum_{j<k}\zeta_{j,i;k,i},
\end{align*}
where we used $w^*(\infty) = 1$ in (\ref{eqn: correction alpha}) and (\ref{eqn: correction beta}).
   Hence, we have the \textit{limiting variability equations}:
 \begin{eqnarray}\label{eqn: variability equations - nonmatrix}
 c^2_{a,i} & = &\sum_{j = 1}^{K} (\lambda_{j,i}/\lambda_i) c^2_{a,j,i} + (\lambda_{0,i}/\lambda_i) c^2_{a,0,i} + c^2_{\beta_i},  \nonumber \\
 c^2_{a,i,j}& = & p_{i,j}c^2_{d,i} + (1-p_{i,j}) + c^2_{\alpha_{i,j}}, \nonumber \\
 c^2_{d,i} & = & c^2_{a,i}, \quad 1 \le i \le K.
\end{eqnarray}
where we used the fact that $w_{i}(t) \to 1$ as $t \to \infty$.

For a concise matrix notation, let
\begin{align*}
\mathbf{I}(t) & \equiv (I_{a,1}(t),\dots,I_{a,K}(t), I_{a,1,1}(t),\dots, I_{a,K,K}(t), I_{d,1}(t),\dots,I_{d,K}(t)), \\
\mathbf{b}(t) & \equiv (b_{a,1}(t),\dots,b_{a,K}(t), b_{a,1,1}(t),\dots, b_{a,K,K}(t), b_{d,1}(t),\dots,b_{d,K}(t)), \\
\mathbf{M}(t) & \equiv (M_{m,n}(t))\in \RR^{(2K + K^2)^2},\quad m,n \in\{a_1,\dots,a_K,a_{1,1},\dots, a_{K,K}, d_1,\dots,d_K\}, \\
\mathbf{c}^2 &  \equiv(c^2_{a,1},\dots, c^2_{a,K},c^2_{a,1,1}, \dots, c^2_{a,K,K}, c^2_{d,i},\dots, c^2_{d,K}),
\end{align*}
where
\begin{align*}
& b_{a,i}(t)\equiv \frac{\lambda_{0,i}}{\lambda_i} I_{a,0,i}(t) + \beta_i(t), \quad b_{a,i,j} \equiv (1-p_{i,j}) + \alpha_{i,j}(t),  \\
 & b_{d,i}(t) \equiv (1-w_i(t)) I_{s,i}(t); \quad  M_{a_i, a_{j,i}(t)} = \frac{\lambda_{j,i}}{\lambda_i}, \\
 & M_{a_{i,j},d_i}(t) = p_{i,j}, M_{d_i, a_i}(t) = w_{i}(t),\qandq M_{m,n}(t) = 0 \hbox{ otherwise}.
\end{align*}
Then the IDC equations can be expressed concisely as
\begin{equation}\label{eqn: IDC equations}
(\mathbf{E} - \mathbf{M}(t))\mathbf{I}(t) = \mathbf{b}(t),
\end{equation}
while the limiting variability equations can be expressed as
\begin{equation}\label{eqn: limiting variability equations}
(\mathbf{E} - \mathbf{M}(\infty))\mathbf{c}^2 = \mathbf{b}(\infty),
\end{equation}
where $\mathbf{E}\in \RR^{(2K + K^2)^2}$ is the identity matrix.

The following theorem states that these equations have unique solutions.
\begin{theorem}\label{Thm: solution to IDC equations}
Assume that $I-P'$ is invertible. Then $\mathbf{E} - \mathbf{M}(t)$ is invertible for each fixed $t\in\RR^+\cup \{\infty\}$.  Hence, for any given $t$ and $\mathbf{b}$, the IDC equations in {\em \eqn{eqn: IDC equations}} have the unique solution $$\mathbf{I}(t) = (\mathbf{E} - \mathbf{M}(t))^{-1}\mathbf{b}(t)$$ and the limiting variability equations in
{\em \eqn{eqn: limiting variability equations}} have the unique solution $$\mathbf{c} = (\mathbf{E} - \mathbf{M}(\infty))^{-1}\mathbf{b}(\infty).$$
\end{theorem}
 \paragraph{Proof.} Let $\delta_{i,j}$ be the Kronecker delta function.  Then substituting the equations for $I_{a,j,i}(t)$ and $I_{d,i}(t)$ into the equation for $I_{a,i}(t)$, we obtain an equation set for $I_{a,i}(t)$ with coefficient matrix
$\left(\delta_{i,j} -(\lambda_{j,i}/\lambda_i)p_{j,i}w_{j}(t)\right) \in \RR^{K^2}$.  Note that $(\lambda_{j,i}/\lambda_i)w_{j}(t) \le 1$ for  $t\in\RR^+\cup \{\infty\}$, the invertibility of $I-P'$ implies that the equations for $I_{a,i}(t)$ have an unique solution.  Substituting in the solution for $I_{a,i}(t)$, we obtain solutions for $I_{a,i,j}(t)$ and $I_{d,i}(t)$.
$\square$


\begin{remark}\label{Rmk: connect Kim}
{\em
(The Kim \cite{K11b,K11a} MMPP$(2)$ decomposition.)
In Kim \cite{K11b,K11a}, a
decomposition approximation of queueing networks based on MMPP(2)/GI/1 queues was investigated. MMPP(2) stands for Markov modulated Poission process with 2 underlying states.
The four rate parameters in the MMPP(2) are determined from the approximations of the mean, IDC and the third moment process of the arrival process at a pre-selected time $t_0$ and the limiting variability parameter of the arrival process.
The IDC and third moment processes are approximated by the network equations with correction terms motivated from the Markovian routing settings.

At first glance, the IDC equations proposed here are quite similar to the network equations used in \cite{K11b}, see (20), (22) and (31) there.
However, our method are different in three aspects.
First, our approach does not fit the flows to special processes (MMPP in \cite{K11b}), instead we partially characterize the flows by the IDC and apply the RQ algorithm reviewed in \S \ref{sec:RQ_review}.
Second, the entire IDC function is utilized in the RQ algorithm, whereas \cite{K11b} used IDC evaluated at a pre-selected time $t_0$ to fit the parameters of the MMPP.
Third, we rely on more detailed heavy-traffic limit to propose asymptotically exact correction terms, see \S 5.3 of the
appendix. 
\hfill $\square$
}
\end{remark}

\subsection{RQNA for Tree-Structured Queueing Networks}\label{sec: RQNA tree}

With the IDC equations developed in \S \ref{sec: IDC equations}, we immediately obtain an elementary algorithm for tree-structured OQNs.
A \textit{tree-structured queueing network} is an OQN whose topology forms a directed tree.  Recall that a directed tree is a connected directed graph whose underlying undirected graph is a tree.  The tree-structured network is a special case of feed-forward network in which the superposed flows at each node have no common origin.

This special structure greatly simplifies the IDC-based RQNA algorithm.
First, there is no customer feedback, which significantly simplify the IDC equations as well as the dependence in the queueing network.
Second, for any internal flow $A_{i,j}$ that is non-zero, we must have $\alpha_{i,j} = 0$ for the correction term in (\ref{eqn: splitting equation}), see discussions in \S 5.3 of the appendix. 
Finally, the tree structure implies that $\beta_{i} = 0$ for the correction term for superposition because all superposed processes are independent.

We summarize the procedure in Algorithm \ref{alg: tree RQNA}.
\begin{algorithm}[h!]
    \SetKwInOut{Require}{Require}
    \SetKwInOut{Output}{Output}
    \Require{The queueing network has tree structure.}
    \Output{Solution to the IDC equations (\ref{eqn: IDC equations}).}

    \For{$i=1$ to $n$}
     {

     	$\lambda_i \leftarrow \lambda_{0,i} + \sum_{j < i} \lambda_{j} p_{j,i}$;

 		$\rho_{i} \leftarrow \lambda_i / \mu_i$;
       	
        $c^2_{a,i}
   \leftarrow \sum_{j < i} \frac{\lambda_{j,i}}{\lambda_i} c^2_{a,j,i} + \frac{\lambda_{0,i}}{\lambda_i} c^2_{a,0,i}$;

  		$c_{x,i}^2 \leftarrow c_{a,i}^2 + c_{s,i}^2$;
  		
  		$w_i(t) \leftarrow w^{*} ((1-\rho_i)^{2} \lambda_i t/(\rho_i c_{x,i}^2))$;

      $I_{a_i}(t) \leftarrow \sum_{j < i}\frac{\lambda_{j,i}}{\lambda_i} \left(p_{j,i}\left(w_j(t) I_{a,j} (t) + (1-w_j(t)) I_{s,j}(t)\right) + (1-p_{j,i})\right) + \frac{\lambda_{0,i}}{\lambda_i} I_{a,0,i}(t)$;

       $I_{d_i}(t) \leftarrow  w_i(t) I_{a,i} (t) + (1-w_i(t)) I_{s,i}(t)$;

        \For{$j<i$}
        {
       		 $I_{a,i,j}(t) \leftarrow p_{i,j}I_{d,i}(t) + (1-p_{i,j})$;
        }

     }

   \Return $\mathbf{I}(t)$.
    \caption{The RQNA algorithm for approximating the IDC's at each time $t$ in a tree-structured queueing network.}\label{alg: tree RQNA}
\end{algorithm}
To elaborate,
with these simplifications of the correction terms, the equations in (\ref{eqn: IDC equations - nonmatrix}),
yield, for $1 \le i, j \le K$,
\begin{align*}
 I_{a_i}(t) & = \sum_{j = 1}^{K} \frac{\lambda_{j,i}}{\lambda_i} I_{a_{j,i}}(t) + (\lambda_{0,i}/\lambda_i) I_{a_{0,i}}(t),   \\
 I_{a_{i,j}}(t)& = p_{i,j}I_{d_i}(t) + (1-p_{i,j}),   \\
 I_{d_i}(t) & = w_i(t) I_{a_i} (t) + (1-w_i(t)) I_{s_i}(t).
\end{align*}

The IDC equations in this setting inherit a special structure that allows a recursive algorithm.  Note that the stations in the tree-structured network can be partitioned into disjoint layers $\{\mathcal{L}_1, \dots, \mathcal{L}_l\}$ such that for station $i \in \mathcal{L}_k$, it takes only the input flows from $j \in \bigcup_{j = 1}^{k-1}\mathcal{L}_j$ for $1\le k \le l$.  To simplify the notation, we sort the node in the order of their layers and assign arbitrary order to nodes within the same layer.  If $i \in \mathcal{L}_k$, then $\bigcup_{j = 1}^{k-1}\mathcal{L}_{j} \subset \{1,2, \dots, i-1\}$, so that $\lambda_{j,i} = 0$ for all $j \ge i$.
Hence, by substituting in the equations for $I_{d_i}$ and $I_{a_{i,j}}$ into that of $I_{a_i}$, we have
\begin{align}
I_{a_i}(t)
  & = \sum_{j = 1}^{K} \frac{\lambda_{j,i}}{\lambda_i} \left(p_{j,i}\left(w_j(t) I_{a_j} (t) + (1-w_j(t)) I_{s_j}(t)\right) + (1-p_{j,i})\right) + \frac{\lambda_{0,i}}{\lambda_i} I_{a_{0,i}}(t),\notag \\
  & = \sum_{j < i}\frac{\lambda_{j,i}}{\lambda_i} \left(p_{j,i}\left(w_j(t) I_{a_j} (t) + (1-w_j(t)) I_{s_j}(t)\right) + (1-p_{j,i})\right) + \frac{\lambda_{0,i}}{\lambda_i} I_{a_{0,i}}(t).\label{dep approx in feed-forward}
\end{align}
Note that (\ref{dep approx in feed-forward}) exhibits a lower-triangular shape so that we can explicitly write down the solution in the order of the stations.

\section{Feedback Elimination}\label{sec: feedback elimination}
In this section, we discuss the case in which customers can return (feedback) to a queue after receiving service there.
Customer feedback introduces dependence between the arrival process and the service times, even when the service times themselves are mutually independent.  As a result, the decomposition $I_w(t) = I_a(t) + c_s^2$ in \eqn{eqn: IDW decomposition} is no longer valid.
Indeed, assuming that it is, as we have done so far, can introduce serious errors, as we show in our simulation examples.
We address this problem by introducing a feedback elimination procedure.  We start with the so-called immediate feedback in \S \ref{sec: immediate} and generalize it into near-immediate feedback in \S\ref{sec: near immediate}.

\subsection{Immediate Feedback Elimination}\label{sec: immediate}
In Section III of \cite{WW83qna} it is observed that it is often helpful to pre-process the model data by eliminating immediate feedback for queues with feedback.  We now show how that can be done for the RQNA algorithm.

We consider a single queue with i.i.d. feedback. In this case, all feedback is \textit{immediate feedback}, meaning that the customer feeds back to the same queue immediately after completing service, without first going through another service station.  For a $GI/GI/1$ model allowing feedback, all
feedback is necessarily immediate because there is
only one queue.

Normally, the immediate feedback returns the customer back to the end of the queue.  However, in the immediate feedback elimination procedure, the approximation step is to put the customer back at the head of the line so that the customer receives a geometrically random number of service times all at once.  Clearly this does not alter the queue length process or the workload process, because the approximation step is work-conserving.

The modified system is a single-server queue with a new service-time distribution and without feedback.  Let $N_p$ denote a geometric random variable with success probability $1-p$ and support $\mathbb{N}^+$, the positive natural numbers, then the new service time can be expressed as
\begin{equation}
 S_p = \sum_{i = 1}^{N_p} S_i,
\end{equation}
where $S_i$'s are i.i.d. copies of the original service times.
 This modification in service times results in a change in the service scv.  By the conditional variance formula, the scv of the total service time is $\tilde{c}_s^2 = p + (1-p) c_s^2$.  The new service IDC in the modified system is the IDC of the stationary renewal process associated with the new service times.  To obtain the new service IDC, we need only find the Laplace Transform of the new service distribution, then apply the algorithm in \S 3.1 of the appendix.
 We provides the details in \S 4 of the appendix.

For the mean waiting time, we need to adjust for per-visit waiting time by multiplying the waiting time in the modified system by $(1-p)$.  Note that $(1-p)^{-1}$ is the mean number of visits by a customer in the original system.

In \S 4.1 of \cite{WY20flow} it is shown that the modified system after the immediate feedback elimination procedure shares the same HT limits of the queue length process, the external departure process, the workload process and the waiting time process.
Hence, the immediate feedback elimination procedure as an approximation is asymptotically exact in the heavy-traffic limit.

\subsection{Near-Immediate Feedback}\label{sec: near immediate}

Now, we consider general OQNs, where the feedback does not necessarily happen immediately, meaning that a departing customer
 may visit other queues before coming back to the feedback queue. To treat general OQNs,
 we extend the immediate feedback concept to the \textit{near-immediate feedback}, which depends on the traffic intensities of the queues on the path the customer took before feedback happens.
 The near-immediate feedback is defined as any feedback that does not go through any queue with higher traffic intensity.

By default, the RQNA algorithm eliminates all near-immediate feedback. To help understand near-immediate feedback, consider a modified OQN with one bottleneck queue, denoted by $h$.  A \textit{bottleneck queue} is a queue with the highest traffic intensity in the network. While all non-bottleneck queues have service times set to $0$ so that they serve as instantaneous switches.  In the reduced network, we define an external arrival $\hat{A}_0$ to the bottleneck queue to be any external arrival that arrive at the bottleneck queue for the first time.  Hence, an external arrival may have visited one or multiple non-bottleneck queues before its first visit to the bottleneck queue. In particular, the external arrival process can be expressed as the superposition of (i) the original external arrival process $A_{0,h}$ at station $h$; and (ii) the Markov splitting of the external arrival process $A_{0,i}$ at station $i$ with probability $\hat{p}_{i,h}$, for $i\neq h$, where $\hat{p}_{i,h}$ denote the probability of a customer that enters the original system at station $i$ ends up visiting the bottleneck station $h$.  For the explicit formula of $\hat{p}_{i,h}$, see Remark 3.2 of \cite{WY20flow}.

In \S 4.2 of \cite{WY20flow}, we showed that this reduced network is asymptotically equivalent in the HT limit to the single-server queue with i.i.d. feedback that we considered in \S \ref{sec: immediate}.  In particular, the arrival process of the equivalent single-station system is $\hat{A}_0$ as described above, the service times remain unchanged and the feedback probability is $\hat{p}$, which is exactly the probability of a near-immediate feedback in the original system; see (3.9) of \cite{WY20flow} for the expression of $\hat{p}$.
Hence we showed that eliminating all feedback at the bottleneck queue as described above prior to analysis is asymptotically correct in HT for OQNs with a single bottleneck queue in terms of the queue length process, the external departure process, the workload process and the waiting time process.
Moreover, the different variants of the algorithm  - eliminating all near immediate feedback or only the near-immediate feedback at the bottleneck queues - are asymptotically exact in the HT limit for an OQN with a single-bottleneck queue, because only the bottleneck queues have nondegenerate HT limit.
In contrast, if there are multiple bottleneck queues, the HT limit requires multidimensional RBM, which is not used in our RQNA.

\section{The Full RQNA Algorithm}\label{sec: RQNA alg}

As basic input parameters, the RQNA algorithm requires the model data specified in \S \ref{secData}:
\begin{enumerate}
  \item Network topology specified by the routing matrix $P$;
  \item External arrival processes specified by (i) the interarrival distribution, if renewal; or (ii) rate $\lambda$ and IDC; or (iii) a realized sample path of the stationary external arrival process;
  \item Service renewal process specified by (i) the service distribution; or (ii) the rate and IDC; or (iii) a realized sample path of the stationary service renewal process.
\end{enumerate}

Combining the traffic-rate equation, the limiting variability equation, the IDC equation and the feedback elimination procedure, we have obtained a general framework for the RQNA algorithm, which we summarize in Algorithm \ref{alg: general RQNA}.
We remark that the RQNA algorithm becomes much simpler in the case without customer feedbacks, see \S \ref{sec: RQNA tree} for more discussion.

\begin{algorithm}
    \SetKwInOut{Require}{Require}
    \SetKwInOut{Output}{Output}
    \Require{Specification of the correction terms  $\alpha_{i,j}(t)$ in \S \ref{sec: splitting}
    and $\beta_i(t)$ in \S
    \ref{sec: superposition}, a set of stations to perform feedback elimination as specified in \S \ref{sec: feedback elimination}
    and the flows to eliminate for each of the selected station.}
    \Output{Approximation of the system performance measures.}

	Solve the traffic rate equations by $\lambda = (I - P')^{-1}\lambda_0$ as in \S \ref{secRateEqns}
	and let $\rho_i = \lambda_i/\mu_i$; \label{alg RQNA step rate equation}

    Solve the limiting variability equations by $\mathbf{c} = (\mathbf{E} - \mathbf{M}(\infty))^{-1}\mathbf{b}(\infty)$
    specified in \S \ref{sec: IDC equations}; \label{alg RQNA step var par}

    Solve the IDC equations by $\mathbf{I}(t) = (\mathbf{E} - \mathbf{M}(t))^{-1}\mathbf{b}(t)$ for the total arrival IDCs,
    where we use $\mathbf{c}$ from Step \ref{alg RQNA step var par} in (\ref{depIDC2}); \label{alg RQNA step IDC equation}

    Select a set of stations to perform feedback elimination, as in \S \ref{sec: feedback elimination}.  For each selected station, identify the flows to eliminate, then identify the corresponding feedback probability, the modified service IDC as in \S \ref{sec: immediate} as well as the reduced network.  Repeat Step \ref{alg RQNA step rate equation} to Step \ref{alg RQNA step IDC equation} on the reduced network to obtain the modified IDW (as the sum of the modified total arrival IDC and the modified service scv) at the selected station.

    Apply the RQ algorithm in \eqn{eqn: RQ approx} to
    obtain the approximations for the mean steady-state workload at each station. \label{alg step RQ}

    Apply the formulas in Remark \ref{Rmk:other_perf} and \ref{rmNet} to obtain approximations for the expected values of the steady-state queue length and  waiting time
    at each queue and the total sojourn time for the system. \label{alg step others}

    \caption{A general framework of the RQNA algorithm for the approximation of the system performance measures.}\label{alg: general RQNA}
\end{algorithm}

The general framework here allows different choices of
(i) the correction terms $\alpha_{i,j}$ in \S \ref{sec: splitting} and $\beta_i$ in \S
    \ref{sec: superposition}
    and (ii) the feedback elimination procedure.
The default correction terms are given in (\ref{eqn: correction alpha}) and (\ref{eqn: correction beta}).
For the feedback elimination procedure, we apply near-immediate feedback elimination to all stations.
In \S 6 of the appendix 
we discuss an additional tuning function to fine tune the performance of our RQNA algorithm.


\section{Numerical Studies}\label{sec: numerical}
In this section, we discuss examples of networks with significant near-immediate feedback from \cite{DNR94}.
We show that the near-immediate feedback in these examples makes a big difference in the performance descriptions.
Hence our predictions with and without feedback elimination are very different.
We find that our RQNA with near-immediate feedback elimination performs as well or better than the other algorithms.
Additional numerical examples appear in our previous papers and in \S 7 of the appendix. 

\subsection{A Three-Station Example}
\begin{figure}[h!]
  \centering
  \begin{tikzpicture}[auto, >=latex']
    \everymath{\displaystyle}
    \tikzstyle{every picture}+=[remember picture]

    \tikzstyle{coor-format} = [coordinate]
    \tikzstyle{b-format} = [draw, thick,shape=datastore,inner sep=0.1cm,minimum width=1.5cm]
    \tikzstyle{s1-format} = [rectangle, draw=black, fill=black!60]
    \tikzstyle{s2-format} = [rectangle, draw=black, fill=black!30]
    \tikzstyle{p-format} = [circle, draw=black, thick, minimum height=0.6cm, minimum width=0.6cm]

    \node [coor-format, label=90:{$\lambda_{0,1} = 0.225$}] (In1) {};
    \node [coor-format,  right of=In1, node distance=1.5cm] (Split1) {};
    \node [coor-format,  right of=In1, node distance=1cm] (Mid1) {};
    \node [b-format, right of=Split1, node distance=0.8cm] (Buf1) {Queue 1};
    \draw [draw,->, thick] (In1) -- (Buf1);

    \node [p-format, right of=Buf1, node distance=1.6cm] (Pool1) {};
    \draw [draw,->, thick] (Buf1) -- (Pool1);

    \node [coor-format, right of=Pool1, node distance=0.5cm,
    label=140:{}] (Mid2) {};

        \node [b-format, right of=Mid2, node distance=1.4cm] (Buf2) {Queue 2};

    \draw [draw, ->, thick] (Pool1) -- (Buf2);

    \node [coor-format, below of = Mid1, label=45:{}](low1){};

    \path [draw, ->, thick] (low1) -- (Mid1);

    \node [p-format, right of=Buf2, node distance=1.6cm] (Pool2) {};
    \draw [draw, ->, thick] (Buf2) -- (Pool2);

    \node [coor-format, right of=Pool2, node distance=0.5cm,
    label=75:{$p_{2,3}$}] (Mid3) {};

    \node [b-format, right of = Mid3, node distance=1.75cm,
    label=140:{}] (Buf3) {Queue 3};
    \draw [draw,->, thick] (Pool2) -- (Buf3);

    \node [coor-format, below of = Mid3, label=45:{$p_{2,1}$}](low2){};

    \path [draw, -, thick] (Mid3) |- (low2) |- (low1);

    \node [p-format, right of=Buf3, node distance=1.6cm] (Pool3) {};
    \draw [draw,->, thick] (Buf3) -- (Pool3);

    \node [coor-format, right of=Pool3, node distance=1.4cm,
    label=140:{$p_{3,2}$}] (Out3) {};
    \draw [draw,->, thick] (Pool3) -- (Out3);

    \node [coor-format, right of=Pool3, node distance=0.5cm,
    label=140:{ }] (Mid4) {};

    \node [coor-format, above of=Mid4, node distance=1cm,
    label=140:{ }] (up1) {};

    \node [coor-format, above of=Mid2, node distance=1cm,
    label=140:{ }] (up2) {};

    \path [draw, -, thick] (Mid4) |- (up1) |- (up2);
    \draw [draw,->, thick] (up2) -- (Mid2);

 \end{tikzpicture}
    \caption{A three-station example.}\label{fig: 3 stations}
\end{figure}
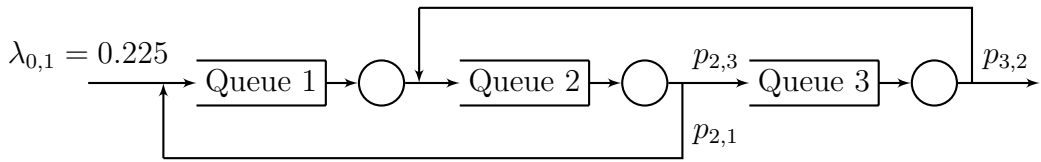

In this section, we look at the suite of three-station examples \S 3.1 of \cite{DNR94} depicted in Figure \ref{fig: 3 stations}.  This example is designed to have three stations that are tightly coupled with each other, so that the dependence among the queues and the flows is fairly complicated.

In this example, we have three stations in tandem but also allow customer feedback from station 2 to station 1 and from station 3 to station 2, with probability $p_{2,1} = p_{2,3} = p_{3,2} = 0.5$.  The only external arrival process is a Poisson process which arrives at station 1 with rate $\lambda_{0,1} = 0.225$, hence by (\ref{traffic_rate_equations}) the effective arrival rate is $\lambda_1 = 0.675, \lambda_2 = 0.9$ and $\lambda_3 = 0.45$.

For the service distributions, we consider the same sets of parameters as in \cite{DNR94}, summarized in Table \ref{table: rho} and \ref{table: scv}.  Note that Case 2 is relatively more challenging because there are two bottlneck stations; in contrast, all the other cases have only one.

\begin{table}[h!]
\parbox{.45\linewidth}{
\caption
{Traffic intensity of the four cases in the three-station example.}\label{table: rho}
\vspace{10pt}
\centering
\begin{tabular}{c|ccc}
\hline
Case & $\rho_1$ & $\rho_2$ & $\rho_3$ \\
\hline
1 & 0.675 & 0.900 & 0.450 \\
2 & 0.900 & 0.675 & 0.900 \\
3 & 0.900 & 0.675 & 0.450 \\
4 & 0.900 & 0.675 & 0.675 \\
\hline
\end{tabular}}
\hspace{30pt}
\parbox{.45\linewidth}{
\caption
{Variability of the service distributions of the four cases in the three-station example.}\label{table: scv}
\centering
\begin{tabular}{c|ccc}
\hline
Case & $c_{s,1}^2$ & $c_{s,2}^2$ & $c_{s,3}^2$ \\
\hline
A & 0.00 & 0.00 & 0.00 \\
B & 2.25 & 0.00 & 0.25 \\
C & 0.25 & 0.25 & 2.25 \\
D & 0.00 & 2.25 & 2.25 \\
E & 8.00 & 8.00 & 0.25 \\
\hline
\end{tabular}}
\end{table}

We now compare the RQNA approximations and four previous algorithms as in \S 7.3 of the appendix,
with the simulated mean sojourn times at each station, as well as total sojourn time of the network.   The sojourn time for each station is defined as the waiting time plus the service time at that station, whereas the total sojourn time of the network is defined as in (\ref{eqn: total sojourn}).  We consider two cases of the RQNA algorithm: (1) the plain RQNA algorithm without feedback elimination, as in Algorithm \ref{alg: general RQNA} and (2) the RQNA algorithm with feedback elimination, as discussed in \S \ref{sec: feedback elimination}.

For RQNA with feedback elimination, we apply feedback elimination to each station that has at least one feedback flow that only passes through stations with equal or lower traffic intensities.  We eliminate all such flows in the feedback elimination procedure.  Take Case 1 for example, we do not apply feedback elimination for Station 1 because all feedback customers go through Station 2, which has higher traffic intensity; we will, however, eliminate the flow from 2 to 1 as well as the flow from 3 to 2 for Station 2, since both Station 1 and 3 have lower traffic intensities.  As another example, for both Station 2 and 3 in case 4, we eliminate the flow from 3 to 2, but we do not eliminate the flow from 2 to 1, since Station 2 and 3 share the same traffic intensity while Station 1 has higher traffic intensity.

Tables \ref{table: three stations} and \ref{table: three stations case D} expand Tables II and III in \cite{DNR94} by adding values for
 (1) the mean total sojourn time and (2) the RQ and RQNA approximations, with and without feedback elimination.  For each table, we indicate by an asterisk in the last column the stations where elimination is applied.

We observed that the plain RQNA algorithm works well for stations with moderate to low traffic intensities, but not so satisfactory for congested stations.  On the other hand, the accuracy of the RQNA algorithm with feedback elimination is on par with, if not better than the best previous algorithm.

\begin{table}
\scriptsize
	\caption
{A comparison of six approximation methods to simulation for the total sojourn time in the three-station example in Figure \ref{fig: 3 stations} with parameters specified in Table \ref{table: rho} and \ref{table: scv}.\label{table: three stations}}
\centering
\begin{tabular}{ll|llll|lll}
\hline
\multicolumn{2}{c}{Case} & Simulation & QNA & QNET & SBD & RQ & RQNA & RQNA (elim) \\
\hline
A   & 1 & 40.39 (3.75\%)  & 20.5 (-49\%)  & diverging     & 43.0 (6.4\%)   & 73.9 (83\%)  & 83.5 (107\%)  & 44.8 (11.0\%) \\
    & 2 & 59.58 (3.29\%)  & 36.0 (-40\%)  & 56.7 (-4.9\%) & 58.2 (-2.4\%)  & 78.0 (31\%)  & 94.3 (58\%)   & 69.3 (16.4\%) \\
    & 3 & 40.72 (4.78\%)  & 24.0 (-41\%)  & 38.7 (-5.0\%) & 40.2 (-1.3\%)  & 57.2 (41\%)  & 74.7 (83\%)   & 43.3 (6.3\%)  \\
    & 4 & 42.12 (3.36\%)  & 26.2 (-38\%)  & 41.8 (-0.7\%) & 42.7 (1.3\%)   & 59.3 (41\%)  & 75.1 (78\%)   & 41.2 (-2.2\%) \\
\hline
B   & 1 & 52.40 (2.64\%)  & 42.0 (-20\%)  & 52.6 (0.4\%)  & 50.2 (-4.2\%)  & 72.4 (38\%)  & 93.7 (79\%)  & 53.1 (1.4\%)  \\
    & 2 & 91.52 (3.77\%)  & 94.1 (2.8\%)  & 83.7 (-8.5\%) & 95.3 (4.1\%)   & 109  (20\%)  & 169 (85\%)  & 94.5 (3.2\%)  \\
    & 3 & 61.68 (3.44\%)  & 72.2 (17\%)   & 61.9 (0.4\%)  & 60.9 (-1.3\%)  & 79.4 (29\%)  & 133 (115\%) & 60.5 (-1.9\%) \\
    & 4 & 63.34 (2.83\%)  & 75.8 (20\%)   & 64.1 (1.3\%)  & 64.7 (2.1\%)   & 83.0 (31\%)  & 135 (113\%) & 62.4 (-1.4\%) \\
\hline
C   & 1 & 44.24 (1.96\%)  & 31.3 (-29\%)  & 37.0 (-16\%)  & 47.1 (6.4\%)   & 75.7 (71\%)  & 91.4 (106\%)  & 42.1 (-4.8\%) \\
    & 2 & 92.42 (4.23\%)  & 87.4 (-5.4\%) & 91.2 (-1.4\%) & 91.6 (-0.83\%) & 106  (15\%)  & 156 (68\%)      & 96.0 (3.8\%)  \\
    & 3 & 44.26 (4.69\%)  & 33.2 (-25\%)  & 44.0 (-0.7\%) & 45.0 (1.7\%)   & 61.3 (38\%)  & 84.2 (90\%)     & 44.0 (-0.6\%) \\
    & 4 & 50.20 (1.04\%)  & 41.4 (-18\%)  & 51.1 (1.7\%)  & 52.2 (4.0\%)   & 67.4 (34\%)  & 91.2 (82\%)     & 45.9 (-8.6\%) \\
\hline
E   & 1 & 134.4 (4.77\%)  & 265 (97\%)    & 155 (15\%)    & 116 (-14\%)    & 158 (17\%)   & 305 (127\%)     & 120 (-11\%)   \\
    & 2 & 213.1 (3.47\%)  & 308 (45\%)    & 228 (7.1\%)   & 206 (-3.3\%)   & 234 (10\%)   & 367 (72\%)      & 173 (-19\%)   \\
    & 3 & 138.7 (3.97\%)  & 244 (76\%)    & 161 (16\%)    & 135 (-2.5\%)   & 163 (17\%)   & 300 (116\%)     & 136 (-2.0\%)  \\
    & 4 & 155.1 (4.37\%)  & 252 (63\%)    & 168 (8.2\%)   & 147 (-5.0\%)   & 178 (15\%)   & 312 (101\%)     & 148 (-4.8\%)  \\
\hline
\end{tabular}

\end{table}

\begin{table}
\scriptsize
	\caption
{A comparison of six approximation methods to simulation for the sojourn time at each station of the three-station example in Figure \ref{fig: 3 stations} for Case D in Table \ref{table: rho} and \ref{table: scv}.}\label{table: three stations case D}
\centering
\begin{tabular}{ll|llll|lll}
 \hline
 Case & Station & Simulation & QNA & QNET & SBD & RQ & RQNA & RQNA (elim) \\
 \hline
 D1  & 1     & 2.476 (0.61\%) & 2.24 (-9.4\%)  & 2.48 (0.3\%)  & 2.47 (-0.1\%) & 2.47 (-0.28\%) & 2.68 (7.8\%)   & 2.68 (7.8\%)\\
     & 2     & 10.85 (3.21\%) & 14.9 (37\%)    & 11.6 (6.5\%)   & 11.4 (5.2\%)   & 19.8 (83\%)    & 28.4 (162\%)   & 11.1$^*$ (2.7\%) \\
     & 3     & 2.544 (0.63\%) & 2.53 (-0.8\%) & 2.54 (-0.0\%) & 2.59 (1.6\%)   & 2.57 (1.2\%)   & 2.53 (-0.7\%) & 2.53 (-0.7\%) \\
     & Total & 55.81 (2.58\%) & 71.4 (28\%)    & 58.8 (5.3\%)   & 58.2 (4.3\%)   & 91.8 (64\%)    & 127 (127\%)    & 57.6 (3.3\%) \\
 \hline
 D2  & 1     & 11.35 (3.29\%) & 8.01 (-29\%)   & 10.8 (-4.5\%)  & 11.1 (-1.9\%)  & 13.7 (20\%)    & 16.6 (46\%)    & 11.3$^*$ (0.1\%)\\
     & 2     & 2.643 (1.25\%) & 2.96 (12\%)    & 2.75 (4.0\%)   & 2.82 (6.7\%)   & 2.85 (7.8\%)   & 3.06 (16\%)    & 3.06 (16\%) \\
     & 3     & 26.87 (2.04\%) & 32.9 (22\%)    & 26.8 (-0.4\%) & 24.9 (-7.5\%)  & 27.5 (2.2\%)   & 36.4 (35\%)    & 31.1$^*$ (16\%) \\
     & Total & 98.36 (1.82\%) & 102  (3.4\%)   & 97.2 (-1.2\%)  & 94.4 (-4.0\%)  & 104  (6.0\%)   & 132 (34\%)     & 105 (7.1\%)\\
 \hline
 D3  & 1     & 11.39 (3.04\%) & 7.95 (-30\%)   & 11.0 (-3.5\%)  & 11.3 (-0.5\%) & 15.8 (39\%)    & 16.5 (45\%)    & 11.3$^*$ (-0.5\%)\\
     & 2     & 2.290 (1.27\%) & 2.90 (27\%)    & 2.53 (10\%)    & 2.26 (-1.4\%)  & 2.57 (12\%)    & 3.04 (33\%)    & 2.10$^*$ (-8.2\%) \\
     & 3     & 2.220 (0.59\%) & 2.40 (7.9\%)   & 2.38 (7.0\%)   & 2.59 (16\%)    & 2.39 (7.6\%)   & 2.43 (9.6\%)   & 2.43 (9.6\%) \\
     & Total & 47.72 (2.51\%) & 40.2 (-16\%)   & 47.8 (0.2\%)  & 48.2 (1.0\%)   & 62.6 (31\%)    & 66.6 (39\%)    & 47.5 (0.51\%) \\
  \hline
 D4  & 1     & 11.30 (6.39\%) & 7.97 (-29\%)   & 10.9 (-3.2\%)  & 11.3 (0.3\%)  & 14.2 (26\%)    & 16.43 (45\%)   & 11.3$^*$ (0.3\%)\\
     & 2     & 2.414 (1.12\%) & 2.93 (21\%)    & 2.64 (9.5\%)   & 2.60 (7.7\%)   & 2.65 (10\%)    & 3.05 (26\%)    & 2.10$^*$ (-13\%) \\
     & 3     & 5.886 (1.05\%) & 6.83 (16\%)    & 6.31 (7.3\%)   & 6.17 (4.8\%)   & 6.47 (10\%)    & 6.85 (16\%)    & 5.95$^*$ (1.1\%) \\
     & Total & 55.24 (4.37\%) & 49.3 (-11\%)   & 56.0 (1.4\%)   & 56.7 (2.7\%)   & 69.3 (25\%)    & 75.5 (37\%)    & 54.3 (-1.7\%)  \\
 \hline
    \multicolumn{3}{c}{Average absolute relative error}  & 20.24\% & 4.72\% & 4.52\% & 21.61\% & 42.60\% & 5.51\% \\
 \hline
\end{tabular}
\end{table}

\subsection{A Ten-Station Example}
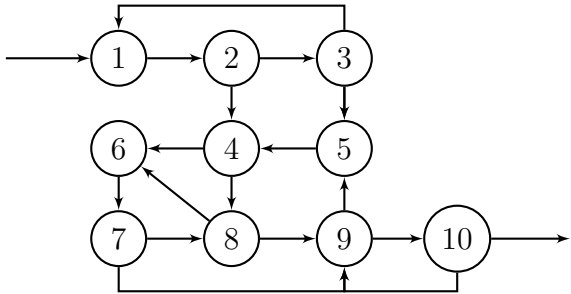
\begin{figure}[h!]
  \centering
  \begin{tikzpicture}[auto, >=latex']
    \everymath{\displaystyle}
    \tikzstyle{every picture}+=[remember picture]

    \tikzstyle{coor-format} = [coordinate]
    \tikzstyle{b-format} = [draw, thick,shape=datastore,inner sep=0.1cm,minimum width=1.5cm]
    \tikzstyle{s1-format} = [rectangle, draw=black, fill=black!60]
    \tikzstyle{s2-format} = [rectangle, draw=black, fill=black!30]
    \tikzstyle{p-format} = [circle, draw=black, thick, minimum height=0.6cm, minimum width=0.6cm]

    \node [coor-format] (In) {};

    \node [p-format, right of=In1, node distance=1.5cm] (Q1) {1};

    \draw [draw,->, thick] (In1) -- (Q1);

	\node [p-format, right of=Q1, node distance=1.5cm] (Q2) {2};
	\node [p-format, right of=Q2, node distance=1.5cm] (Q3) {3};	
	\node [p-format, below of=Q1, node distance=1.2cm] (Q6) {6};	
	\node [p-format, right of=Q6, node distance=1.5cm] (Q4) {4};	
	\node [p-format, right of=Q4, node distance=1.5cm] (Q5) {5};	
	\node [p-format, below of=Q6, node distance=1.2cm] (Q7) {7};	
	\node [p-format, right of=Q7, node distance=1.5cm] (Q8) {8};	
	\node [p-format, right of=Q8, node distance=1.5cm] (Q9) {9};	
	\node [p-format, right of=Q9, node distance=1.5cm] (Q10) {10};
	
	\node [coor-format, right of=Q10, node distance=1.5cm] (Out) {};
	
	\node [coor-format, above of=Q1, node distance=0.7cm] (Cor1) {};
	\node [coor-format, above of=Q3, node distance=0.7cm] (Cor2) {};
	\path [draw, ->, thick] (Q3) |- (Cor2) |- (Cor1) -> (Q1);	
	
	\draw [draw,->, thick] (Q1) -- (Q2);	
	\draw [draw,->, thick] (Q2) -- (Q3);	
	\draw [draw,->, thick] (Q2) -- (Q4);	
	\draw [draw,->, thick] (Q3) -- (Q5);	
	\draw [draw,->, thick] (Q3) -- (Q5);	
	\draw [draw,->, thick] (Q5) -- (Q4);	
	\draw [draw,->, thick] (Q4) -- (Q6);	
	\draw [draw,->, thick] (Q4) -- (Q8);	
	\draw [draw,->, thick] (Q6) -- (Q7);	
	\draw [draw,->, thick] (Q7) -- (Q8);	
	\draw [draw,->, thick] (Q8) -- (Q6);	
	\draw [draw,->, thick] (Q8) -- (Q9);	
	\draw [draw,->, thick] (Q9) -- (Q5);	
	\draw [draw,->, thick] (Q9) -- (Q10);	
	\draw [draw,->, thick] (Q10) -- (Out);	
	
	\node [coor-format, below of=Q7, node distance=0.7cm] (Cor3) {};	
	\node [coor-format, below of=Q9, node distance=0.7cm] (Cor4) {};
	\node [coor-format, below of=Q10, node distance=0.7cm] (Cor5) {};

	\path [draw, ->, thick] (Q10) |- (Cor5) |- (Cor4) -> (Q9);	
	\path [draw, ->, thick] (Q7) |- (Cor3) |- (Cor4) -> (Q9);	

	\end{tikzpicture}
	 \caption{A ten-station with customer feedback example.}\label{fig: 10 stations}
\end{figure}

We conclude with the 10-station OQN example with feedback considered in \S 3.5 of \cite{DNR94}.
It is depicted here in Figure \ref{fig: 10 stations}.

The only exogenous arrival process is Poission with rate 1.  For each station, if there are two routing destinations, the departing customer follows Markovian routing with equal probability, each being $0.5$.  The vector of mean service times is $(0.45, 0.30, 0.90, 0.30, 0.38571,$ $0.20, 0.1333, 0.20, 0.15, 0.20)$, so that the traffic intensity vector is $(0.6, 0.4, 0.6, 0.9, 0.9,0.6,$ $0.4, 0.6, 0.6,  0.4)$. The scv's at these stations are $(0.5, 2, 2, 0.25, 0.25, 2, 1, 2, 0.5, 0.5)$, where we assume a Erlang distribution if $c_s^2 < 1$, an exponential distribution if $c_s^2 = 1$ and a hyperexponential distribution if $c_s^2 > 1$.

In particular, note that stations 4 and 5 are bottleneck queues, having equal traffic intensity,
far greater than the traffic intensities at the other queues.  Moreover, these two stations are quite closely coupled.
Thus, at first glance, we expect that SBD with two-dimensional RBM should perform very well, which proves to be correct.
Moreover, this example should be challenging for RQNA because it is based on heavy-traffic limits for OQNs with only a single bottleneck, thus involving only one-dimensional RBM.

In Table \ref{table: 10 stations}, we report the simulation estimates and approximattions for the steady-state mean sojourn time (waiting time plus service time) at each station, as well as the total sojourn time of the system, calculated as in \eqn{eqn: total sojourn}.  For the approximations, we compare QNA from \cite{WW83qna}, QNET from \cite{HN90}, SBD from \cite{DNR94}, RQ from \cite{WY18RQ1} (with estimated IDC), as well as the RQNA algorithms here.  The simulation, QNA, QNET and SBD columns are taken from Table XIV of \cite{DNR94}.

Again, we consider two versions of RQNA algorithm, the first one does not eliminate feedback, while the second one (marked by `elim') applies the feedback elimination procedure.  As before, in eliminating customer feedback, for each station, we identify the near-immediate feedback flows as the flows that come back to the station after completing service, without passing through any station with a higher traffic intensity.  We then eliminate all near-immediate feedback flows, apply plain RQNA algorithm on the reduced network and use the new RQNA approximation as the approximation for that station.

\begin{table}
\scriptsize
\caption
{A comparison of six approximation methods to simulation for the mean steady-state sojourn times at each station of the open queueing network in Figure \ref{fig: 10 stations}.}\label{table: 10 stations}
\centering
\begin{tabular}{l|llll|lll}
\hline
Station & Simulation & QNA & QNET & SBD & RQ & RQNA & RQNA (elim) \\
 \hline
1     & 0.99 (0.86\%)  & 0.97 (-2.8\%)  & 1.00 (0.2\%)   & 1.00  (0.4\%) & 0.97 (-2.0\%)  & 1.09 (9.2\%)   &  1.00$^*$ (0.4\%) \\
2     & 0.55 (0.69\%)  & 0.58 (6.0\%)   & 0.56 (2.6\%)   & 0.55 (0.2\%)  & 0.55 (-0.1\%)  & 0.56 (1.3\%)   &  0.56 (1.4\%) \\
3     & 2.82 (1.93\%)  & 2.93 (4.2\%)   & 2.90 (3.2\%)   & 2.76 (-2.0\%) & 2.96 (5.0\%)   & 3.40 (21\%)    &  2.75$^*$ (-2.5\%) \\
4     & 1.79 (3.71\%)  & 1.34 (-25\%)   & 1.41 (-21\%)   & 1.76 (-1.6\%) & 2.34 (31\%)    & 3.51 (97\%)    &  2.11$^*$ (18\%) \\
5     & 2.92 (4.77\%)  & 2.49 (-15\%)   & 2.44 (-17\%)   & 2.81 (-3.6\%) & 3.77 (29\%)    & 9.07 (211\%)   &  3.35$^*$ (15\%) \\
6     & 0.58 (0.78\%)  & 0.64 (10\%)    & 0.62 (7.4\%)   & 0.59 (2.2\%)  & 0.60 (3.8\%)   & 0.70 (20\%)  &  0.49$^*$ (-16\%) \\
7     & 0.24 (0.28\%)  & 0.24 (-1.7\%)  & 0.26 (7.1\%)   & 0.27 (11\%)   & 0.23 (-3.0\%)  & 0.24 (-1.3\%) &  0.24 (-1.3\%) \\
8     & 0.58 (0.67\%)  & 0.64 (9.6\%)   & 0.61 (4.6\%)   & 0.60 (1.7\%)  & 0.61 (3.9\%)   & 0.70 (20\%)  &  0.59$^*$ (0.6\%) \\
9     & 0.34 (0.63\%)  & 0.32 (-6.1\%)  & 0.35 (2.0\%)   & 0.43 (26\%)   & 0.33 (-4.2\%)  & 0.73 (111\%)   &  0.42$^*$ (21\%) \\
10    & 0.29 (0.19\%)  & 0.30 (2.4\%)   & 0.29 (1.4\%)   & 0.28 (-1.7\%) & 0.28 (-1.5\%)  & 0.26 (-8.7\%)  &  0.26 (-8.7\%) \\
Total & 22.0 (2.45\%)  & 20.3 (-7.9\%)  & 20.4 (-7.3\%)  & 22.4 (1.7\%)  & 26.1 (18\%)    & 44.5 (102\%)   &  24.2$^*$ (9.9\%)  \\
\hline
\end{tabular}
\end{table}

We make the following observations from this numerical example:
\begin{enumerate}
\item Particular attention should be given to the two bottleneck stations: 4 and 5.  Note that
QNA and QNET produce $15-25\%$ error, which is satisfactory, but SBD does far better with only $1-4\%$ error.
  \item The RQNA algorithm without feedback elimination can perform very poorly with high traffic intensity and high feedback probability, presumably due to the break down of the IDW decomposition in (\ref{eqn: IDW decomposition}).
  \item With feedback elimination, the RQNA algorithm performs significantly better and is competitive with previous algorithms in this complex setting, producing $15-18\%$ error at stations 4 and 5.  The performance of RQNA at the tightly coupled bottleneck queues evidently suffers because the current RQNA depends heavily on one-dimensional RBM.
\end{enumerate}

\section{Conclusions}\label{secConclusions}

In this paper we developed a new robust queueing network queueing analyzer (RQNA) based on indices of dispersion.
The indices of dispersion are scaled variance time curves.
They enable the approximations to exploit the dependence in the arrival processes over time to describe
the mean workload as a function of the traffic intensity at each station.

After reviewing the indices of dispersion and
the robust queueing approximation for a single queue in \S \ref{sec:review_IDC_RQ}, we developed the important
variability linear equations for the IDCs of the internal arrival processes in
\S \ref{sec: IDC framework}.  We then introduced the extra step of feedback elimination in \S \ref{sec: feedback elimination}.
These approximations draw heavily on heavy-traffic limits in \cite{WY18RQ1,WY18Dep,WY20flow} involving one-dimensional RBM.
We put all this together into a full algorithm in \S \ref{sec: RQNA alg}, developing a simplified version for
networks with tree structure in \S\ref{sec: RQNA tree}.

We then evaluated the performance of the new RQNA-IDC by making comparisons with simulations
for various examples in \S \ref{sec: numerical} and \S 7 of the appendix.
These experiments confirm that
RQNA-IDC is remarkably effective.  There are many excellent directions for future research, including (i) developing
refined approximations for the flows that exploit multi-dimensional RBM instead of just one-dimensional RBM
and (ii) extending RQNA-IDC to other OQN models, e.g., with multiple servers and other service disciplines.

\section*{Acknowledgements}
This work was done while the Wei You was a graduate student at Columbia University,
where both authors received support from NSF grant CMMI 1634133.

\bibliographystyle{informs2014}
\bibliography{refs}

\end{document}